\documentclass[12pt,a4paper]{amsart}

\usepackage{amsmath,amsthm,amssymb}
\usepackage{graphicx}
\usepackage{tikz-cd}
\usepackage{paralist}
\usepackage{environ}
\usepackage{xcolor}
\usepackage{centernot}
\usepackage{mathtools}
\usepackage[utf8]{inputenc}
\usepackage[cm]{fullpage}
\usepackage{graphics}
\usepackage{amscd}
\usepackage{mathrsfs}
\usepackage{amsfonts}
\usepackage{lscape}
\usepackage{tikz}
\usetikzlibrary{matrix,calc}
\usepackage{multirow}
\usepackage{color}
\usepackage{graphicx}
\usepackage{marvosym}
\usepackage{mdframed}
\usepackage{paralist}
\usepackage{multicol}

\usepackage{pstricks,pst-text,pst-grad,pst-node,pst-3dplot,pstricks-add,pst-poly,pst-coil} 
\usepackage{pst-fun,pst-blur} 

\usepackage{hyperref}
\newtheorem{theorem}{Theorem}
\theoremstyle{plain}

\newtheorem{definition}{Definition}
\newtheorem{example}{Example}

\newtheorem{proposition}{Proposition}

\begin{document}

\title[Determinant of cubic-matrix of orders 2 and 3]{The Determinant of cubic-matrix of order 2 and order 3, some basic properties and algorithms}

\author[Armend Salihu]{Armend Salihu}
\address{Armend Salihu: Department of Computer Science, Faculty of Contemporary Sciences, South East European University, Tetovo, North Macedonia}
\email{ar.salihu@gmail.com}

\author[Orgest ZAKA]{Orgest ZAKA}
\address{Orgest ZAKA: Department of Mathematics-Informatics, Faculty of Economy and Agribusiness, Agricultural University of Tirana, Tirana, Albania}
\email{ozaka@ubt.edu.al, gertizaka@yahoo.com}

%\dedicatory{ }


\begin{abstract}
Based on geometric intuition, in this paper we are trying to give an idea and visualize the meaning of the determinants for the cubic-matrix. 
In this paper we have analysed the possibilities of developing the concept of determinant of matrices with three indexed 3D Matrices.
We define the concept of determinant for cubic-matrix of order 2 and order 3, study and prove some basic properties for calculations of determinants of cubic-matrix of order 2 and 3.
Furthermore we have also tested several square determinant properties and noted that these properties also are applicable on this concept of 3D Determinants. 
\end{abstract}


\maketitle
\tableofcontents

\section{Introduction and Preliminaries}
Based on the determinant of 2D square matrices \cite{Salihu6, Salihu7, Salihu8, ArtinM, BretscherO, Schneide-Barker, Lang}, as well as determinant of rectangular matrices \cite{Salihu1, Salihu2, Salihu3, Salihu4, Salihu5, Amiri-etal, Radic1, Radic2, MAKAREWICZetal} we have come to the idea of developing the concept of determinant of 3D cubic matrices, our concept is based on permutation expansion method. Encouraged by geometric intuition, in this paper we are trying to give an idea and visualize the meaning of the determinants for the cubic-matrix. 
Our early research mainly lies between geometry, algebra, matrix theory, etc., (see \cite{PetersZakaDyckAM}, \cite{ZakaDilauto}, \cite{ZakaFilipi2016}, \cite{FilipiZakaJusufi}, \cite{ZakaCollineations}, \cite{ZakaVertex}, \cite{ZakaThesisPhd}, \cite{ZakaPetersIso}, \cite{ZakaPetersOrder}, \cite{ZakaMohammedSF}, \cite{ZakaMohammedEndo}).

This paper is continuation of the ideas that arise based on previous researches of 3D matrix ring with element from any whatever field $F$ see \cite{zaka3DmatrixRing}, but here we study the case when the field $F$ is the field of real numbers $\mathbb{R}$. In this paper we follow a different method from the calculation of determinants of 3D matrix, which is studied in \cite{zaka3DGLnnp}. In contrast to the meaning of the determinant as a multi-scalar studied in \cite{zaka3DGLnnp}, in this paper we give a new definition, for the determinant of the 3D-cubic-matrix, which is a real-number.

In the papers \cite{zaka3DmatrixRing, zaka3DGLnnp}, have been studied in detail, properties for 3D-matrix, therefore, those studied properties are also valid for 3D-cubic-Matrix.

Our point in this paper is to provide a concept of determinant of 3D matrices. Our concept is based on Milne-Thomson \cite{Milne-Thomson} or permutation method used in regular square matrices.

The following is definition of 3D matrices provided by Zaka in 2017 (see \cite{zaka3DmatrixRing, zaka3DGLnnp}):
\begin{definition}
3-dimensional $m\times n\times p$ matrix will call, a matrix which has: \emph{m-horizontal layers} (analogous to m-rows), \emph{n-vertical page} (analogue with n-
columns in the usual matrices) and \emph{p-vertical layers} (p-1 of which are hidden).

The set of these matrixes the write how:
\begin{equation}
M_{m \times n \times p} ({\rm F})=\{a_{i,j,k}|a_{i,j,k} \in F-{\rm field \hspace{6pt} } \forall i=\overline{1,m};\hspace{6pt} j=\overline{1,n};\hspace{6pt} k=\overline{1,p}\}
\end{equation}
\end{definition}

\begin{figure}[h]
\includegraphics[width=0.4\textwidth]{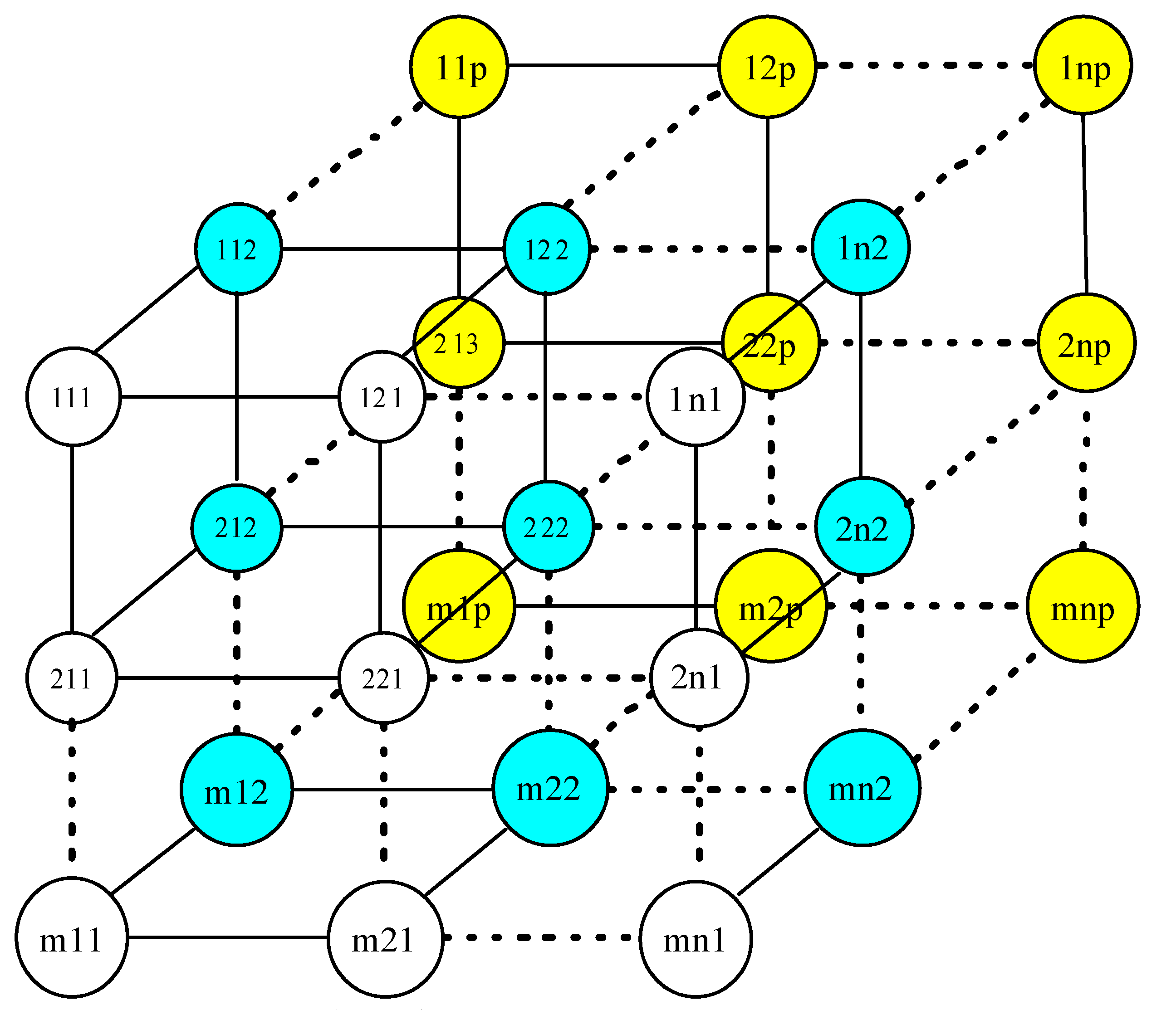}
\centering
\caption{3D-Matrix view}
\end{figure}

In the following is presented the determinant of 3D-cubic matrices, as well several properties which are adopted from 2D square determinants. 

\section{Cubic-Matrix of Order 2 and 3 and their Determinants} %$ $\\

A cubic-matrix $A_{n \times n \times n}$ for $n=2,3, \ldots$, called "cubic-matrix of order $n$". 

For $n=1$ we have that the cubic-matrix of order 1 is an element of $F$.

Let us now consider the set of cubic-matrix of order $n$, for $n=2$ or $n=3$, with elements from a field $F$ (so when cubic-matrix of order $n$, there are: $n-$vertical pages, $n-$horizontal layers and $n-$vertical layers).

From \cite{zaka3DmatrixRing, zaka3DGLnnp} we have that, the addition of 3D-matrix stands also for cubic-matrix of order 2 and 3. Also, the set of cubic-matrix of order 2 and 3 forms an commutative group (Abelian Group) related to 3Dmatrix addition.

\subsection{Determinants of Cubic-Matrix of Order 2 and 3} $ $\\

In this section, we will define and describe the meaning of the determinants of cubic-matrix of order 2 and order 3, with elements from a field $F$. Recall that a cubic-matrix $A_{n \times n \times n}$ for $n=2,3, \ldots$, called "cubic-matrix of order $n $". 

For $n=1$ we have that the cubic-matrix of order 1 is an element of $F$.

Let us now consider the set of cubic-matrix of order $n$, with elements from a field $F$ (so when cubic-matrix of order $n$, there are: $n-$vertical pages, $n-$horizontal layers and $n-$vertical layers),

\[
\mathcal{M}_n(F)=\{ A_{n \times n \times n}=(a_{ijk})_{n \times n \times n} | a_{ijk} \in F, \forall i=\bar{1,n}; j=\bar{1,n}; k=\bar{1,n} \}
\]

In this paper, we define the \emph{determinant of cubic-matrix} as a element from this field, so the map, 

\[
\begin{aligned}
\det : \mathcal{M}_n(F) & \rightarrow F \\
\forall A \in \mathcal{M}_n(F) & \mapsto \det(A) \in F
\end{aligned}
\]

Below we give two definitions, how we will calculate the determinant of the cubic-matrix of order 2 and order 3.

\begin{definition}\label{defOrder2}
Let $A \in \mathcal{M}_2(F)$ be a $2 \times 2 \times 2$, with elements from a field $F$. 
\[
A_{2 \times 2\times 2}=
\begin{pmatrix}
\left.\begin{matrix}
a_{111} & a_{121}\\ 
a_{211} & a_{221}
\end{matrix}\right|
\begin{matrix}
a_{112} & a_{122}\\ 
a_{212} & a_{222}
\end{matrix}
\end{pmatrix}
\]

Determinant of this cubic-matrix, we called,
\[
\det[A_{2 \times 2\times 2}]=\det
\begin{pmatrix}
\left.\begin{matrix}
a_{111} & a_{121}\\ 
a_{211} & a_{221}
\end{matrix}\right|
\begin{matrix}
a_{112} & a_{122}\\ 
a_{212} & a_{222}
\end{matrix}
\end{pmatrix}=a_{111}\cdot a_{222} - a_{112}\cdot a_{221} - a_{121}\cdot a_{212} + a_{122}\cdot a_{211}
\]
\end{definition}

The follow example is case where cubic-matrix, is with elements from the number field $\mathbb{R}$.
\begin{example}
Let's have the cubic-matrix, with element in number field $\mathbb{R}$, 
\[
A_{2 \times 2\times 2}=
\begin{pmatrix}
\left.\begin{matrix}
2 & 1\\ 
3 & 5
\end{matrix}\right|
\begin{matrix}
4 & 7\\ 
3 & 2
\end{matrix}
\end{pmatrix}
\]

then according to the definition \ref{defOrder2}, we calculate the Determinant of this cubic-matrix, and have,
\[
\det[A_{2 \times 2\times 2}]=\det
\begin{pmatrix}
\left.\begin{matrix}
2 & 1\\ 
3 & 5
\end{matrix}\right|
\begin{matrix}
4 & 7\\ 
3 & 2
\end{matrix}
\end{pmatrix}=2\cdot 2 - 4\cdot 5 - 1\cdot 3 + 7\cdot 3=2.
\]

\end{example}

We are trying to expand the meaning of the determinant of cubic-matrix, for order 3 (so when cubic-matrix, there are: 3-vertical pages, 3-horizontal layers and 3-vertical layers).

\begin{definition} \label{defOrder3}
Let $A \in \mathcal{M}_3(F)$ be a $3 \times 3 \times 3$ cubic-matrix with element from a field $F$, 
\[
A_{3 \times 3\times 3}=
\begin{pmatrix} \left.
\begin{matrix}
a_{111} & a_{121} & a_{131}\\ 
a_{211} & a_{221} & a_{231}\\ 
a_{311} & a_{321} & a_{331}
\end{matrix}\right|
\begin{matrix}
a_{112} & a_{122} & a_{132}\\ 
a_{212} & a_{222} & a_{232}\\ 
a_{312} & a_{322} & a_{332}
\end{matrix}\left|
\begin{matrix}
a_{113} & a_{123} & a_{133}\\ 
a_{213} & a_{223} & a_{233}\\ 
a_{313} & a_{323} & a_{333}
\end{matrix}\right.
\end{pmatrix} 
. \]

Determinant of this cubic-matrix, we called,

\begin{equation}
\det[A_{3 \times 3\times 3}]=
\det\begin{pmatrix} \left.
\begin{matrix}
a_{111} & a_{121} & a_{131}\\ 
a_{211} & a_{221} & a_{231}\\ 
a_{311} & a_{321} & a_{331}
\end{matrix}\right|
\begin{matrix}
a_{112} & a_{122} & a_{132}\\ 
a_{212} & a_{222} & a_{232}\\ 
a_{312} & a_{322} & a_{332}
\end{matrix}\left|
\begin{matrix}
a_{113} & a_{123} & a_{133}\\ 
a_{213} & a_{223} & a_{233}\\ 
a_{313} & a_{323} & a_{333}
\end{matrix}\right.
\end{pmatrix} 
\end{equation}
$$=a_{111} \cdot a_{222} \cdot a_{333} - a_{111} \cdot a_{232} \cdot a_{323} - a_{111} \cdot a_{223} \cdot a_{332} + a_{111} \cdot a_{233} \cdot a_{322}$$
$$-a_{112} \cdot a_{221} \cdot a_{333} + a_{112} \cdot a_{223} \cdot a_{331} + a_{112} \cdot a_{231} \cdot a_{323} - a_{112} \cdot a_{233} \cdot a_{321}$$
$$+a_{113} \cdot a_{221} \cdot a_{332} - a_{113} \cdot a_{222} \cdot a_{331} - a_{113} \cdot a_{231} \cdot a_{322} + a_{113} \cdot a_{232} \cdot a_{321}$$
$$-a_{121} \cdot a_{212} \cdot a_{333} + a_{121} \cdot a_{213} \cdot a_{332} + a_{121} \cdot a_{232} \cdot a_{313} - a_{121} \cdot a_{233} \cdot a_{312}$$
$$+a_{122} \cdot a_{211} \cdot a_{333} - a_{122} \cdot a_{213} \cdot a_{331} - a_{122} \cdot a_{231} \cdot a_{313} + a_{122} \cdot a_{233} \cdot a_{311}$$
$$-a_{123} \cdot a_{211} \cdot a_{332} + a_{123} \cdot a_{212} \cdot a_{331} + a_{123} \cdot a_{231} \cdot a_{312} - a_{123} \cdot a_{232} \cdot a_{311}$$
$$+a_{131} \cdot a_{212} \cdot a_{323} - a_{131} \cdot a_{213} \cdot a_{322} - a_{131} \cdot a_{222} \cdot a_{313} + a_{131} \cdot a_{223} \cdot a_{312}$$
$$-a_{132} \cdot a_{211} \cdot a_{323} + a_{132} \cdot a_{213} \cdot a_{321} + a_{132} \cdot a_{221} \cdot a_{313} - a_{132} \cdot a_{223} \cdot a_{311}$$
$$+a_{133} \cdot a_{211} \cdot a_{322} - a_{133} \cdot a_{212} \cdot a_{321} - a_{133} \cdot a_{221} \cdot a_{312} + a_{133} \cdot a_{222} \cdot a_{311}$$

\end{definition}

The follow example is case where cubic-matrix, is with elements from the number field $\mathbb{R}$.
\begin{example}
Let's have the cubic-matrix of order 3, with element from number field (field of real numbers) $\mathbb{R}$,
\[
A_{3 \times 3\times 3}=
\det\begin{pmatrix}
\left.
\begin{matrix}
1&4&2\\2&0&0\\0&4&2
\end{matrix}\right |
\begin{matrix}
3&1&3\\5&1&3\\3&2&0
\end{matrix}\left |
\begin{matrix}
2&1&0\\0&1&0\\2&1&0
\end{matrix} \right.
\end{pmatrix}
. \]

Then, we calculate the Determinant of this cubic-matrix following the Definition \ref{defOrder3}, and have that,

\begin{equation*}
\det[A_{3 \times 3\times 3}]=
\det\begin{pmatrix}
\left.
\begin{matrix}
1&4&2\\2&0&0\\0&4&2
\end{matrix}\right |
\begin{matrix}
3&1&3\\5&1&3\\3&2&0
\end{matrix}\left |
\begin{matrix}
2&1&0\\0&1&0\\2&1&0
\end{matrix} \right.
\end{pmatrix}
\end{equation*}
$$= 1 \cdot 1 \cdot 0 - 1 \cdot 1 \cdot 0 - 1 \cdot 3 \cdot 1 + 1 \cdot 0 \cdot 2 - 3 \cdot 0 \cdot 0 + 3 \cdot 1 \cdot 2 + 3 \cdot 0 \cdot 1 - 3 \cdot 0 \cdot 4 + 2 \cdot 0 \cdot 0- 2 \cdot 1 \cdot 2- 2 \cdot 0 \cdot 2 +2 \cdot 3 \cdot 4 $$
$$ - 4 \cdot 5 \cdot 0+4 \cdot 0 \cdot 0+4 \cdot 3 \cdot 2-4 \cdot 0 \cdot 3+1 \cdot 2 \cdot 0-1 \cdot 0 \cdot 2-1 \cdot 0 \cdot 2+1 \cdot 0 \cdot 0-1 \cdot 2 \cdot 0+1 \cdot 5 \cdot 2+1 \cdot 0 \cdot 3 -1 \cdot 3 \cdot 0$$
$$+2 \cdot 5 \cdot 1-2 \cdot 0 \cdot 2-2 \cdot 1 \cdot 2+2 \cdot 1 \cdot 3-3 \cdot 2 \cdot 1+3 \cdot 0 \cdot 4+3 \cdot 0 \cdot 2-3 \cdot 1 \cdot 0+0 \cdot 2 \cdot 2- 0\cdot 5 \cdot 4-0 \cdot 0 \cdot 3+0 \cdot 1 \cdot 0$$
$$=-3+6-4+24+24+10+10-4+6-6=63$$
\end{example}

%============PROPERTIES of DET=================

\section{Some basic properties of determinants for cubic-matrix of order 2 and order 3} %$ $\\

\begin{definition} We will call $I_n$, a unit-3D-cubic-matrix of order 2 or 3, with elements $e_{ijk}$, which are:
\[
e_{ijk}=\left\{
\begin{matrix}
0 & for, \quad i\neq j\neq k\\
1 & for, \quad i=j=k
\end{matrix}
\right.
\]
\end{definition}

\begin{proposition} 
For every unit-cubic-matrix of order 2 or 3, with element from number field $\mathbb{R}$, we have that $\det(I_n)=1$.
\end{proposition}
\proof
Let's have the unit cubic-matrix of order 2, 
\[
I_2=
\begin{pmatrix} \left.
\begin{matrix}
1&0\\0&0
\end{matrix}\right| 
\begin{matrix}
0&0\\0&1
\end{matrix}
\end{pmatrix}
\]
then, this determinant is,
\[
\det(I_2)=\det
\begin{pmatrix} \left.
\begin{matrix}
1&0\\0&0
\end{matrix}\right| 
\begin{matrix}
0&0\\0&1
\end{matrix}
\end{pmatrix}
= 1 \cdot 1=1.
\]
Now lets have the unit-cubic matrix of order 3, 
\[
I_3=
\begin{pmatrix} \left.
\begin{matrix}
1&0&0\\0&0&0\\0&0&0
\end{matrix}\right| 
\begin{matrix}
0&0&0\\0&1&0\\0&0&0
\end{matrix} \left|
\begin{matrix}
0&0&0\\0&0&0\\0&0&1
\end{matrix} \right.
\end{pmatrix}
\]
then this determinant is,
\[
\det(I_3)=\det[I_{3 \times 3\times 3}]=
\det\begin{pmatrix} \left.
\begin{matrix}
1&0&0\\0&0&0\\0&0&0
\end{matrix}\right| 
\begin{matrix}
0&0&0\\0&1&0\\0&0&0
\end{matrix} \left|
\begin{matrix}
0&0&0\\0&0&0\\0&0&1
\end{matrix} \right.
\end{pmatrix}=1\cdot 1\cdot 1=1
\]
\qed

\begin{theorem} Suppose that $A$ is a cubic-matrix of order 2 or 3, with a plan where every entry is zero then its determinant is 'zero', so $\det(A)=0.$
\end{theorem}

\proof
1. For plan $i=1$:

Let $A$ be cubic-matrix of order 2, where all elements on the plan $i=1$ are equal to zero, then based on definition \ref{defOrder2}:
\[
\det[A_{2 \times 2\times 2}]=\det
\begin{pmatrix}
\left.\begin{matrix}
0 & 0\\ 
a_{211} & a_{221}
\end{matrix}\right|
\begin{matrix}
0 & 0\\ 
a_{212} & a_{222}
\end{matrix}
\end{pmatrix}=0\cdot a_{222} - 0\cdot a_{221} - 0\cdot a_{212} + 0\cdot a_{211}=0.
\]

2. For plan $i=2$:

Let $A$ be cubic-matrix of order 2, where all elements on the plan $i=2$ are equal to zero, then based on definition \ref{defOrder2}:
\[
A_{2 \times 2\times 2}=
\begin{pmatrix}
\left.\begin{matrix}
a_{111} & a_{121}\\ 
0 & 0
\end{matrix}\right|
\begin{matrix}
a_{112} & a_{122}\\ 
0 & 0
\end{matrix}
\end{pmatrix} =a_{111}\cdot0 - a_{112}\cdot 0 - a_{121}\cdot 0 + a_{122}\cdot 0=0.
\]

3. For plan $j=1$:

Let $A$ be cubic-matrix of order 2, where all elements on the plan $j=1$ are equal to zero, then based on definition \ref{defOrder2}:
\[
\det[A_{2 \times 2\times 2}]=\det
\begin{pmatrix}
\left.\begin{matrix}
0 & a_{121}\\ 
0 & a_{221}
\end{matrix}\right|
\begin{matrix}
0 & a_{122}\\ 
0 & a_{222}
\end{matrix}
\end{pmatrix}=0\cdot a_{222} - 0\cdot a_{221} - a_{121}\cdot 0 + a_{122}\cdot 0 = 0.
\]

4. For plan $j=2$:

Let $A$ be cubic-matrix of order 2, where all elements on the plan $j=2$ are equal to zero, then based on definition \ref{defOrder2}:
\[
\det[A_{2 \times 2\times 2}]=\det
\begin{pmatrix}
\left.\begin{matrix}
a_{111} & 0\\ 
a_{211} & 0
\end{matrix}\right|
\begin{matrix}
a_{112} & 0\\ 
a_{212} & 0
\end{matrix}
\end{pmatrix}=a_{111}\cdot 0 - a_{112}\cdot 0 - 0\cdot a_{212} + 0\cdot a_{211} = 0.
\]

5. For plan $k=1$:

Let $A$ be cubic-matrix of order 2, where all elements on the plan $k=1$ are equal to zero, then based on definition \ref{defOrder2}:
\[
\det[A_{2 \times 2\times 2}]=\det
\begin{pmatrix}
\left.\begin{matrix}
0 & 0\\ 
0 & 0
\end{matrix}\right|
\begin{matrix}
a_{112} & a_{122}\\ 
a_{212} & a_{222}
\end{matrix}
\end{pmatrix}=0\cdot a_{222} - a_{112}\cdot 0 - 0\cdot a_{212} + a_{122}\cdot 0 = 0.
\]

6. For plan $k=2$:

Let $A$ be cubic-matrix of order 2, where all elements on the plan $k=2$ are equal to zero, then based on definition \ref{defOrder2}:
\[
\det[A_{2 \times 2\times 2}]=\det
\begin{pmatrix}
\left.\begin{matrix}
a_{111} & a_{121}\\ 
a_{211} & a_{221}
\end{matrix}\right|
\begin{matrix}
0 & 0\\ 
0 & 0
\end{matrix}
\end{pmatrix}=a_{111}\cdot 0 - 0\cdot a_{221} - a_{121}\cdot 0 + 0\cdot a_{211} = 0.
\]

Now we will consider for third order cubic-matrix, as following.

1. For plan $i=1$:

Let $A$ be cubic-matrix of order 3, where all elements on the plan $i=1$ are equal to zero, then based on definition \ref{defOrder3}:
\[
\det[A_{3 \times 3\times 3}]=
\det\begin{pmatrix} \left.
\begin{matrix}
0 & 0 & 0 \\ 
a_{211} & a_{221} & a_{231}\\ 
a_{311} & a_{321} & a_{331}
\end{matrix}\right|
\begin{matrix}
0 & 0 & 0 \\ 
a_{212} & a_{222} & a_{232}\\ 
a_{312} & a_{322} & a_{332}
\end{matrix}\left|
\begin{matrix}
0 & 0 & 0 \\ 
a_{213} & a_{223} & a_{233}\\ 
a_{313} & a_{323} & a_{333}
\end{matrix}\right.
\end{pmatrix} 
\]
$$=0 \cdot a_{222} \cdot a_{333} - 0 \cdot a_{232} \cdot a_{323} - 0 \cdot a_{223} \cdot a_{332} + 0 \cdot a_{233} \cdot a_{322}$$
$$-0 \cdot a_{221} \cdot a_{333} + 0 \cdot a_{223} \cdot a_{331} + 0 \cdot a_{231} \cdot a_{323} - 0 \cdot a_{233} \cdot a_{321}$$
$$+0 \cdot a_{221} \cdot a_{332} - 0 \cdot a_{222} \cdot a_{331} - 0 \cdot a_{231} \cdot a_{322} + 0 \cdot a_{232} \cdot a_{321}$$
$$-0 \cdot a_{212} \cdot a_{333} + 0 \cdot a_{213} \cdot a_{332} + 0 \cdot a_{232} \cdot a_{313} - 0 \cdot a_{233} \cdot a_{312}$$
$$+0 \cdot a_{211} \cdot a_{333} - 0 \cdot a_{213} \cdot a_{331} - 0 \cdot a_{231} \cdot a_{313} + 0 \cdot a_{233} \cdot a_{311}$$
$$-0 \cdot a_{211} \cdot a_{332} + 0 \cdot a_{212} \cdot a_{331} + 0 \cdot a_{231} \cdot a_{312} - 0 \cdot a_{232} \cdot a_{311}$$
$$+0 \cdot a_{212} \cdot a_{323} - 0 \cdot a_{213} \cdot a_{322} - 0 \cdot a_{222} \cdot a_{313} + 0 \cdot a_{223} \cdot a_{312}$$
$$-0 \cdot a_{211} \cdot a_{323} + 0 \cdot a_{213} \cdot a_{321} + 0 \cdot a_{221} \cdot a_{313} - 0 \cdot a_{223} \cdot a_{311}$$
$$+0 \cdot a_{211} \cdot a_{322} - 0 \cdot a_{212} \cdot a_{321} - 0 \cdot a_{221} \cdot a_{312} + 0 \cdot a_{222} \cdot a_{311} = 0.$$

2. For plan $i=2$:

Let $A$ be cubic-matrix of order 3, where all elements on the plan $i=2$ are equal to zero, then based on definition \ref{defOrder3}:
\[
\det[A_{3 \times 3\times 3}]=
\det\begin{pmatrix} \left.
\begin{matrix}
a_{111} & a_{121} & a_{131}\\ 
0 & 0 & 0 \\ 
a_{311} & a_{321} & a_{331}
\end{matrix}\right|
\begin{matrix}
a_{112} & a_{122} & a_{132}\\ 
0 & 0 & 0 \\ 
a_{312} & a_{322} & a_{332}
\end{matrix}\left|
\begin{matrix}
a_{113} & a_{123} & a_{133}\\ 
0 & 0 & 0 \\ 
a_{313} & a_{323} & a_{333}
\end{matrix}\right.
\end{pmatrix} 
\]
$$=a_{111} \cdot 0 \cdot a_{333} - a_{111} \cdot 0 \cdot a_{323} - a_{111} \cdot 0 \cdot a_{332} + a_{111} \cdot 0 \cdot a_{322}$$
$$-a_{112} \cdot 0 \cdot a_{333} + a_{112} \cdot 0 \cdot a_{331} + a_{112} \cdot 0 \cdot a_{323} - a_{112} \cdot 0 \cdot a_{321}$$
$$+a_{113} \cdot 0 \cdot a_{332} - a_{113} \cdot 0 \cdot a_{331} - a_{113} \cdot 0 \cdot a_{322} + a_{113} \cdot 0 \cdot a_{321}$$
$$-a_{121} \cdot 0 \cdot a_{333} + a_{121} \cdot 0 \cdot a_{332} + a_{121} \cdot 0 \cdot a_{313} - a_{121} \cdot 0 \cdot a_{312}$$
$$+a_{122} \cdot 0 \cdot a_{333} - a_{122} \cdot 0 \cdot a_{331} - a_{122} \cdot 0 \cdot a_{313} + a_{122} \cdot 0 \cdot a_{311}$$
$$-a_{123} \cdot 0 \cdot a_{332} + a_{123} \cdot 0 \cdot a_{331} + a_{123} \cdot 0 \cdot a_{312} - a_{123} \cdot 0 \cdot a_{311}$$
$$+a_{131} \cdot 0 \cdot a_{323} - a_{131} \cdot 0 \cdot a_{322} - a_{131} \cdot 0 \cdot a_{313} + a_{131} \cdot 0 \cdot a_{312}$$
$$-a_{132} \cdot 0 \cdot a_{323} + a_{132} \cdot 0 \cdot a_{321} + a_{132} \cdot 0 \cdot a_{313} - a_{132} \cdot 0 \cdot a_{311}$$
$$+a_{133} \cdot 0 \cdot a_{322} - a_{133} \cdot 0 \cdot a_{321} - a_{133} \cdot 0 \cdot a_{312} + a_{133} \cdot 0 \cdot a_{311}=0.$$

3. For plan $i=3$:

Let $A$ be cubic-matrix of order 3, where all elements on the plan $i=3$ are equal to zero, then based on definition \ref{defOrder3}:
\[
\det[A_{3 \times 3\times 3}]=
\det\begin{pmatrix} \left.
\begin{matrix}
a_{111} & a_{121} & a_{131}\\ 
a_{211} & a_{221} & a_{231}\\ 
0 & 0 & 0 
\end{matrix}\right|
\begin{matrix}
a_{112} & a_{122} & a_{132}\\ 
a_{212} & a_{222} & a_{232}\\ 
0 & 0 & 0 
\end{matrix}\left|
\begin{matrix}
a_{113} & a_{123} & a_{133}\\ 
a_{213} & a_{223} & a_{233}\\ 
0 & 0 & 0 
\end{matrix}\right.
\end{pmatrix} 
\]
$$=a_{111} \cdot a_{222} \cdot 0 - a_{111} \cdot a_{232} \cdot 0 - a_{111} \cdot a_{223} \cdot 0 + a_{111} \cdot a_{233} \cdot 0 $$
$$-a_{112} \cdot a_{221} \cdot 0 + a_{112} \cdot a_{223} \cdot 0 + a_{112} \cdot a_{231} \cdot 0 - a_{112} \cdot a_{233} \cdot 0 $$
$$+a_{113} \cdot a_{221} \cdot 0 - a_{113} \cdot a_{222} \cdot 0 - a_{113} \cdot a_{231} \cdot 0 + a_{113} \cdot a_{232} \cdot 0 $$
$$-a_{121} \cdot a_{212} \cdot 0 + a_{121} \cdot a_{213} \cdot 0 + a_{121} \cdot a_{232} \cdot 0 - a_{121} \cdot a_{233} \cdot 0 $$
$$+a_{122} \cdot a_{211} \cdot 0 - a_{122} \cdot a_{213} \cdot 0 - a_{122} \cdot a_{231} \cdot 0 + a_{122} \cdot a_{233} \cdot 0 $$
$$-a_{123} \cdot a_{211} \cdot 0 + a_{123} \cdot a_{212} \cdot 0 + a_{123} \cdot a_{231} \cdot 0 - a_{123} \cdot a_{232} \cdot 0 $$
$$+a_{131} \cdot a_{212} \cdot 0 - a_{131} \cdot a_{213} \cdot 0 - a_{131} \cdot a_{222} \cdot 0 + a_{131} \cdot a_{223} \cdot 0 $$
$$-a_{132} \cdot a_{211} \cdot 0 + a_{132} \cdot a_{213} \cdot 0 + a_{132} \cdot a_{221} \cdot 0 - a_{132} \cdot a_{223} \cdot 0 $$
$$+a_{133} \cdot a_{211} \cdot 0 - a_{133} \cdot a_{212} \cdot 0 - a_{133} \cdot a_{221} \cdot 0 + a_{133} \cdot a_{222} \cdot 0 = 0.$$

4. For plan $j=1$:

Let $A$ be cubic-matrix of order 3, where all elements on the plan $j=1$ are equal to zero, then based on definition \ref{defOrder3}:
\[
\det[A_{3 \times 3\times 3}]=
\det\begin{pmatrix} \left.
\begin{matrix}
0 & a_{121} & a_{131}\\ 
0 & a_{221} & a_{231}\\ 
0 & a_{321} & a_{331}
\end{matrix}\right|
\begin{matrix}
0 & a_{122} & a_{132}\\ 
0 & a_{222} & a_{232}\\ 
0 & a_{322} & a_{332}
\end{matrix}\left|
\begin{matrix}
0 & a_{123} & a_{133}\\ 
0 & a_{223} & a_{233}\\ 
0 & a_{323} & a_{333}
\end{matrix}\right.
\end{pmatrix} 
\]
$$=0 \cdot a_{222} \cdot a_{333} - 0 \cdot a_{232} \cdot a_{323} - 0 \cdot a_{223} \cdot a_{332} + 0 \cdot a_{233} \cdot a_{322}$$
$$-0 \cdot a_{221} \cdot a_{333} + 0 \cdot a_{223} \cdot a_{331} + 0 \cdot a_{231} \cdot a_{323} - 0 \cdot a_{233} \cdot a_{321}$$
$$+0 \cdot a_{221} \cdot a_{332} - 0 \cdot a_{222} \cdot a_{331} - 0 \cdot a_{231} \cdot a_{322} + 0 \cdot a_{232} \cdot a_{321}$$
$$-a_{121} \cdot 0 \cdot a_{333} + a_{121} \cdot 0 \cdot a_{332} + a_{121} \cdot a_{232} \cdot 0 - a_{121} \cdot a_{233} \cdot 0 $$
$$+a_{122} \cdot 0 \cdot a_{333} - a_{122} \cdot 0 \cdot a_{331} - a_{122} \cdot a_{231} \cdot 0 + a_{122} \cdot a_{233} \cdot 0 $$
$$-a_{123} \cdot 0 \cdot a_{332} + a_{123} \cdot 0 \cdot a_{331} + a_{123} \cdot a_{231} \cdot 0 - a_{123} \cdot a_{232} \cdot 0 $$
$$+a_{131} \cdot 0 \cdot a_{323} - a_{131} \cdot 0 \cdot a_{322} - a_{131} \cdot a_{222} \cdot 0 + a_{131} \cdot a_{223} \cdot 0 $$
$$-a_{132} \cdot 0 \cdot a_{323} + a_{132} \cdot 0 \cdot a_{321} + a_{132} \cdot a_{221} \cdot 0 - a_{132} \cdot a_{223} \cdot 0 $$
$$+a_{133} \cdot 0 \cdot a_{322} - a_{133} \cdot 0 \cdot a_{321} - a_{133} \cdot a_{221} \cdot 0 + a_{133} \cdot a_{222} \cdot 0 = 0.$$

5. For plan $j=2$:

Let $A$ be cubic-matrix of order 3, where all elements on the plan $j=2$ are equal to zero, then based on definition \ref{defOrder3}:
\[
\det[A_{3 \times 3\times 3}]=
\det\begin{pmatrix} \left.
\begin{matrix}
a_{111} & 0 & a_{131}\\ 
a_{211} & 0 & a_{231}\\ 
a_{311} & 0 & a_{331}
\end{matrix}\right|
\begin{matrix}
a_{112} & 0 & a_{132}\\ 
a_{212} & 0 & a_{232}\\ 
a_{312} & 0 & a_{332}
\end{matrix}\left|
\begin{matrix}
a_{113} & 0 & a_{133}\\ 
a_{213} & 0 & a_{233}\\ 
a_{313} & 0 & a_{333}
\end{matrix}\right.
\end{pmatrix} 
\]
$$=a_{111} \cdot 0 \cdot a_{333} - a_{111} \cdot a_{232} \cdot 0 - a_{111} \cdot 0 \cdot a_{332} + a_{111} \cdot a_{233} \cdot 0 $$
$$-a_{112} \cdot 0 \cdot a_{333} + a_{112} \cdot 0 \cdot a_{331} + a_{112} \cdot a_{231} \cdot 0 - a_{112} \cdot a_{233} \cdot 0 $$
$$+a_{113} \cdot 0 \cdot a_{332} - a_{113} \cdot 0 \cdot a_{331} - a_{113} \cdot a_{231} \cdot 0 + a_{113} \cdot a_{232} \cdot 0 $$
$$-0 \cdot a_{212} \cdot a_{333} + 0 \cdot a_{213} \cdot a_{332} + 0 \cdot a_{232} \cdot a_{313} - 0 \cdot a_{233} \cdot a_{312}$$
$$+0 \cdot a_{211} \cdot a_{333} - 0 \cdot a_{213} \cdot a_{331} - 0 \cdot a_{231} \cdot a_{313} + 0 \cdot a_{233} \cdot a_{311}$$
$$-0 \cdot a_{211} \cdot a_{332} + 0 \cdot a_{212} \cdot a_{331} + 0 \cdot a_{231} \cdot a_{312} - 0 \cdot a_{232} \cdot a_{311}$$
$$+a_{131} \cdot a_{212} \cdot 0 - a_{131} \cdot a_{213} \cdot 0 - a_{131} \cdot 0 \cdot a_{313} + a_{131} \cdot 0 \cdot a_{312}$$
$$-a_{132} \cdot a_{211} \cdot 0 + a_{132} \cdot a_{213} \cdot 0 + a_{132} \cdot 0 \cdot a_{313} - a_{132} \cdot 0 \cdot a_{311}$$
$$+a_{133} \cdot a_{211} \cdot 0 - a_{133} \cdot a_{212} \cdot 0 - a_{133} \cdot 0 \cdot a_{312} + a_{133} \cdot 0 \cdot a_{311} = 0.$$

6. For plan $j=3$:

Let $A$ be cubic-matrix of order 3, where all elements on the plan $j=3$ are equal to zero, then based on definition \ref{defOrder3}:
\[
\det[A_{3 \times 3\times 3}]=
\det\begin{pmatrix} \left.
\begin{matrix}
a_{111} & a_{121} & 0 \\ 
a_{211} & a_{221} & 0 \\ 
a_{311} & a_{321} & 0 
\end{matrix}\right|
\begin{matrix}
a_{112} & a_{122} & 0 \\ 
a_{212} & a_{222} & 0 \\ 
a_{312} & a_{322} & 0 
\end{matrix}\left|
\begin{matrix}
a_{113} & a_{123} & 0 \\ 
a_{213} & a_{223} & 0 \\ 
a_{313} & a_{323} & 0 
\end{matrix}\right.
\end{pmatrix} 
\]
$$=a_{111} \cdot a_{222} \cdot 0 - a_{111} \cdot 0 \cdot a_{323} - a_{111} \cdot a_{223} \cdot 0 + a_{111} \cdot 0 \cdot a_{322}$$
$$-a_{112} \cdot a_{221} \cdot 0 + a_{112} \cdot a_{223} \cdot 0 + a_{112} \cdot 0 \cdot a_{323} - a_{112} \cdot 0 \cdot a_{321}$$
$$+a_{113} \cdot a_{221} \cdot 0 - a_{113} \cdot a_{222} \cdot 0 - a_{113} \cdot 0 \cdot a_{322} + a_{113} \cdot 0 \cdot a_{321}$$
$$-a_{121} \cdot a_{212} \cdot 0 + a_{121} \cdot a_{213} \cdot 0 + a_{121} \cdot 0 \cdot a_{313} - a_{121} \cdot 0 \cdot a_{312}$$
$$+a_{122} \cdot a_{211} \cdot 0 - a_{122} \cdot a_{213} \cdot 0 - a_{122} \cdot 0 \cdot a_{313} + a_{122} \cdot 0 \cdot a_{311}$$
$$-a_{123} \cdot a_{211} \cdot 0 + a_{123} \cdot a_{212} \cdot 0 + a_{123} \cdot 0 \cdot a_{312} - a_{123} \cdot 0 \cdot a_{311}$$
$$+0 \cdot a_{212} \cdot a_{323} - 0 \cdot a_{213} \cdot a_{322} - 0 \cdot a_{222} \cdot a_{313} + 0 \cdot a_{223} \cdot a_{312}$$
$$-0 \cdot a_{211} \cdot a_{323} + 0 \cdot a_{213} \cdot a_{321} + 0 \cdot a_{221} \cdot a_{313} - 0 \cdot a_{223} \cdot a_{311}$$
$$+0 \cdot a_{211} \cdot a_{322} - 0 \cdot a_{212} \cdot a_{321} - 0 \cdot a_{221} \cdot a_{312} + 0 \cdot a_{222} \cdot a_{311} = 0.$$

7. For plan $k=1$:

Let $A$ be cubic-matrix of order 3, where all elements on the plan $k=1$ are equal to zero, then based on definition \ref{defOrder3}:
\[
\det[A_{3 \times 3\times 3}]=
\det\begin{pmatrix} \left.
\begin{matrix}
0 & 0 & 0 \\ 
0 & 0 & 0 \\ 
0 & 0 & 0 
\end{matrix}\right|
\begin{matrix}
a_{112} & a_{122} & a_{132}\\ 
a_{212} & a_{222} & a_{232}\\ 
a_{312} & a_{322} & a_{332}
\end{matrix}\left|
\begin{matrix}
a_{113} & a_{123} & a_{133}\\ 
a_{213} & a_{223} & a_{233}\\ 
a_{313} & a_{323} & a_{333}
\end{matrix}\right.
\end{pmatrix} 
\]
$$=0 \cdot a_{222} \cdot a_{333} - 0 \cdot a_{232} \cdot a_{323} - 0 \cdot a_{223} \cdot a_{332} + 0 \cdot a_{233} \cdot a_{322}$$
$$-a_{112} \cdot 0 \cdot a_{333} + a_{112} \cdot a_{223} \cdot 0 + a_{112} \cdot 0 \cdot a_{323} - a_{112} \cdot a_{233} \cdot 0 $$
$$+a_{113} \cdot 0 \cdot a_{332} - a_{113} \cdot a_{222} \cdot 0 - a_{113} \cdot 0 \cdot a_{322} + a_{113} \cdot a_{232} \cdot 0 $$
$$-0 \cdot a_{212} \cdot a_{333} + 0 \cdot a_{213} \cdot a_{332} + 0 \cdot a_{232} \cdot a_{313} - 0 \cdot a_{233} \cdot a_{312}$$
$$+a_{122} \cdot 0 \cdot a_{333} - a_{122} \cdot a_{213} \cdot 0 - a_{122} \cdot 0 \cdot a_{313} + a_{122} \cdot a_{233} \cdot 0 $$
$$-a_{123} \cdot 0 \cdot a_{332} + a_{123} \cdot a_{212} \cdot 0 + a_{123} \cdot 0 \cdot a_{312} - a_{123} \cdot a_{232} \cdot 0 $$
$$+0 \cdot a_{212} \cdot a_{323} - 0 \cdot a_{213} \cdot a_{322} - 0 \cdot a_{222} \cdot a_{313} + 0 \cdot a_{223} \cdot a_{312}$$
$$-a_{132} \cdot 0 \cdot a_{323} + a_{132} \cdot a_{213} \cdot 0 + a_{132} \cdot 0 \cdot a_{313} - a_{132} \cdot a_{223} \cdot 0 $$
$$+a_{133} \cdot 0 \cdot a_{322} - a_{133} \cdot a_{212} \cdot 0 - a_{133} \cdot 0 \cdot a_{312} + a_{133} \cdot a_{222} \cdot 0 = 0. $$

8. For plan $k=2$:

Let $A$ be cubic-matrix of order 3, where all elements on the plan $k=2$ are equal to zero, then based on definition \ref{defOrder3}:
\[
\det[A_{3 \times 3\times 3}]=
\det\begin{pmatrix} \left.
\begin{matrix}
a_{111} & a_{121} & a_{131}\\ 
a_{211} & a_{221} & a_{231}\\ 
a_{311} & a_{321} & a_{331}
\end{matrix}\right|
\begin{matrix}
0 & 0 & 0 \\ 
0 & 0 & 0 \\ 
0 & 0 & 0 
\end{matrix}\left|
\begin{matrix}
a_{113} & a_{123} & a_{133}\\ 
a_{213} & a_{223} & a_{233}\\ 
a_{313} & a_{323} & a_{333}
\end{matrix}\right.
\end{pmatrix} 
\]
$$=a_{111} \cdot 0 \cdot a_{333} - a_{111} \cdot 0 \cdot a_{323} - a_{111} \cdot a_{223} \cdot 0 + a_{111} \cdot a_{233} \cdot 0 $$
$$-0 \cdot a_{221} \cdot a_{333} + 0 \cdot a_{223} \cdot a_{331} + 0 \cdot a_{231} \cdot a_{323} - 0 \cdot a_{233} \cdot a_{321}$$
$$+a_{113} \cdot a_{221} \cdot 0 - a_{113} \cdot 0 \cdot a_{331} - a_{113} \cdot a_{231} \cdot 0 + a_{113} \cdot 0 \cdot a_{321}$$
$$-a_{121} \cdot 0 \cdot a_{333} + a_{121} \cdot a_{213} \cdot 0 + a_{121} \cdot 0 \cdot a_{313} - a_{121} \cdot a_{233} \cdot 0 $$
$$+0 \cdot a_{211} \cdot a_{333} - 0 \cdot a_{213} \cdot a_{331} - 0 \cdot a_{231} \cdot a_{313} + 0 \cdot a_{233} \cdot a_{311}$$
$$-a_{123} \cdot a_{211} \cdot 0 + a_{123} \cdot 0 \cdot a_{331} + a_{123} \cdot a_{231} \cdot 0 - a_{123} \cdot 0 \cdot a_{311}$$
$$+a_{131} \cdot 0 \cdot a_{323} - a_{131} \cdot a_{213} \cdot 0 - a_{131} \cdot 0 \cdot a_{313} + a_{131} \cdot a_{223} \cdot 0 $$
$$-0 \cdot a_{211} \cdot a_{323} + 0 \cdot a_{213} \cdot a_{321} + 0 \cdot a_{221} \cdot a_{313} - 0 \cdot a_{223} \cdot a_{311}$$
$$+a_{133} \cdot a_{211} \cdot 0 - a_{133} \cdot 0 \cdot a_{321} - a_{133} \cdot a_{221} \cdot 0 + a_{133} \cdot 0 \cdot a_{311} = 0.$$

9. For plan $k=3$:

Let $A$ be cubic-matrix of order 3, where all elements on the plan $k=3$ are equal to zero, then based on definition \ref{defOrder3}:
\[
\det[A_{3 \times 3\times 3}]=
\det\begin{pmatrix} \left.
\begin{matrix}
a_{111} & a_{121} & a_{131}\\ 
a_{211} & a_{221} & a_{231}\\ 
a_{311} & a_{321} & a_{331}
\end{matrix}\right|
\begin{matrix}
a_{112} & a_{122} & a_{132}\\ 
a_{212} & a_{222} & a_{232}\\ 
a_{312} & a_{322} & a_{332}
\end{matrix}\left|
\begin{matrix}
0 & 0 & 0 \\ 
0 & 0 & 0 \\ 
0 & 0 & 0 
\end{matrix}\right.
\end{pmatrix} 
\]
$$=a_{111} \cdot a_{222} \cdot 0 - a_{111} \cdot a_{232} \cdot 0 - a_{111} \cdot 0 \cdot a_{332} + a_{111} \cdot 0 \cdot a_{322}$$
$$-a_{112} \cdot a_{221} \cdot 0 + a_{112} \cdot 0 \cdot a_{331} + a_{112} \cdot a_{231} \cdot 0 - a_{112} \cdot 0 \cdot a_{321}$$
$$+0 \cdot a_{221} \cdot a_{332} - 0 \cdot a_{222} \cdot a_{331} - 0 \cdot a_{231} \cdot a_{322} + 0 \cdot a_{232} \cdot a_{321}$$
$$-a_{121} \cdot a_{212} \cdot 0 + a_{121} \cdot 0 \cdot a_{332} + a_{121} \cdot a_{232} \cdot 0 - a_{121} \cdot 0 \cdot a_{312}$$
$$+a_{122} \cdot a_{211} \cdot 0 - a_{122} \cdot 0 \cdot a_{331} - a_{122} \cdot a_{231} \cdot 0 + a_{122} \cdot 0 \cdot a_{311}$$
$$-0 \cdot a_{211} \cdot a_{332} + 0 \cdot a_{212} \cdot a_{331} + 0 \cdot a_{231} \cdot a_{312} - 0 \cdot a_{232} \cdot a_{311}$$
$$+a_{131} \cdot a_{212} \cdot 0 - a_{131} \cdot 0 \cdot a_{322} - a_{131} \cdot a_{222} \cdot 0 + a_{131} \cdot 0 \cdot a_{312}$$
$$-a_{132} \cdot a_{211} \cdot 0 + a_{132} \cdot 0 \cdot a_{321} + a_{132} \cdot a_{221} \cdot 0 - a_{132} \cdot 0 \cdot a_{311}$$
$$+0 \cdot a_{211} \cdot a_{322} - 0 \cdot a_{212} \cdot a_{321} - 0 \cdot a_{221} \cdot a_{312} + 0 \cdot a_{222} \cdot a_{311} = 0.$$

Based on definition \ref{defOrder2} and definition \ref{defOrder3}, we can see that each term is multiplied with elements of each plan once, hence the proof can be easily seen, as presented above.

\qed

\begin{proposition} Suppose that $A$ is cubic-matrix of order 2 or 3, and let's be $B$ the cubic-matrix of same order with A, which obtained from $A$ by multiplying any single: "horizontal layer" or "vertical page" or "vertical layer" with scalar $\alpha$. Then $\det(B)=\alpha \cdot \det(A)$, if $\alpha \neq 0$.
\end{proposition}

\proof

Case 1. The cubic-matrix A of order 2, (and B has order 2), we will proof the case 1 for each "horizontal layer", "vertical page" and "vertical layer", as following:

1. For plan $i=1$:

Let $A$ be cubic-matrix of order 2, where all elements on the plan $i=1$ are equal to zero, then based on definition \ref{defOrder2}:
\[
\det[B_{2 \times 2\times 2}]=\det
\begin{pmatrix}
\left.\begin{matrix}
\alpha \cdot a_{111} & \alpha \cdot a_{121}\\ 
a_{211} & a_{221}
\end{matrix}\right|
\begin{matrix}
\alpha \cdot a_{112} & \alpha \cdot a_{122}\\ 
a_{212} & a_{222}
\end{matrix}
\end{pmatrix}\]
\[=\alpha \cdot a_{111}\cdot a_{222} - \alpha \cdot a_{112}\cdot a_{221} - \alpha \cdot a_{121}\cdot a_{212} + \alpha \cdot a_{122}\cdot a_{211} = \alpha \cdot \det[A_{2 \times 2\times 2}] 
\]

2. For plan $i=2$:

Let $A$ be cubic-matrix of order 2, where all elements on the plan $i=2$ are equal to zero, then based on definition \ref{defOrder2}:
\[
\det[B_{2 \times 2\times 2}]=\det
\begin{pmatrix}
\left.\begin{matrix}
a_{111} & a_{121}\\ 
\alpha \cdot a_{211} & \alpha \cdot a_{221}
\end{matrix}\right|
\begin{matrix}
a_{112} & a_{122}\\ 
\alpha \cdot a_{212} & \alpha \cdot a_{222}
\end{matrix}
\end{pmatrix}
\]
\[
=a_{111}\cdot \alpha \cdot a_{222} - a_{112}\cdot \alpha \cdot a_{221} - a_{121}\cdot \alpha \cdot a_{212} + a_{122}\cdot \alpha \cdot a_{211} = \alpha \cdot \det[A_{2 \times 2\times 2}] 
\]

3. For plan $j=1$:

Let $A$ be cubic-matrix of order 2, where all elements on the plan $j=1$ are equal to zero, then based on definition \ref{defOrder2}:
\[
\det[B_{2 \times 2\times 2}]=\det
\begin{pmatrix}
\left.\begin{matrix}
\alpha \cdot a_{111} & a_{121}\\ 
\alpha \cdot a_{211} & a_{221}
\end{matrix}\right|
\begin{matrix}
\alpha \cdot a_{112} & a_{122}\\ 
\alpha \cdot a_{212} & a_{222}
\end{matrix}
\end{pmatrix}
\]
\[
=\alpha \cdot a_{111}\cdot a_{222} - \alpha \cdot a_{112}\cdot a_{221} - a_{121}\cdot \alpha \cdot a_{212} + a_{122}\cdot \alpha \cdot a_{211} = \alpha \cdot \det[A_{2 \times 2\times 2}] 
\]

4. For plan $j=2$:

Let $A$ be cubic-matrix of order 2, where all elements on the plan $j=2$ are equal to zero, then based on definition \ref{defOrder2}:
\[
\det[B_{2 \times 2\times 2}]=\det
\begin{pmatrix}
\left.\begin{matrix}
a_{111} & \alpha \cdot a_{121}\\ 
a_{211} & \alpha \cdot a_{221}
\end{matrix}\right|
\begin{matrix}
a_{112} & \alpha \cdot a_{122}\\ 
a_{212} & \alpha \cdot a_{222}
\end{matrix}
\end{pmatrix}
\]
\[
=a_{111}\cdot \alpha \cdot a_{222} - a_{112}\cdot \alpha \cdot a_{221} - \alpha \cdot a_{121}\cdot a_{212} + \alpha \cdot a_{122}\cdot a_{211} = \alpha \cdot \det[A_{2 \times 2\times 2}] 
\]

5. For plan $k=1$:

Let $A$ be cubic-matrix of order 2, where all elements on the plan $k=1$ are equal to zero, then based on definition \ref{defOrder2}:
\[
\det[B_{2 \times 2\times 2}]=\det
\begin{pmatrix}
\left.\begin{matrix}
\alpha \cdot a_{111} & \alpha \cdot a_{121}\\ 
\alpha \cdot a_{211} & \alpha \cdot a_{221}
\end{matrix}\right|
\begin{matrix}
a_{112} & a_{122}\\ 
a_{212} & a_{222}
\end{matrix}
\end{pmatrix}
\]
\[
=\alpha \cdot a_{a_{111}}\cdot a_{222} - a_{112}\cdot \alpha \cdot a_{221} - \alpha \cdot a_{121}\cdot a_{212} + a_{122}\cdot \alpha \cdot a_{211} = \alpha \cdot \det[A_{2 \times 2\times 2}] 
\]

6. For plan $k=2$:

Let $A$ be cubic-matrix of order 2, where all elements on the plan $k=2$ are equal to zero, then based on definition \ref{defOrder2}:
\[
\det[B_{2 \times 2\times 2}]=\det
\begin{pmatrix}
\left.\begin{matrix}
a_{111} & a_{121}\\ 
a_{211} & a_{221}
\end{matrix}\right|
\begin{matrix}
\alpha \cdot a_{112} & \alpha \cdot a_{122}\\ 
\alpha \cdot a_{212} & \alpha \cdot a_{222}
\end{matrix}
\end{pmatrix}
\]
\[
=a_{111}\cdot \alpha \cdot a_{222} - \alpha \cdot a_{112}\cdot a_{221} - a_{121}\cdot \alpha \cdot a_{212} + \alpha \cdot a_{122}\cdot a_{211} = \alpha \cdot \det[A_{2 \times 2\times 2}] 
\]

Case 2. The cubic-matrix A of order 3, (and B has order 3)

1. For plan $i=1$:

Let $A$ be cubic-matrix of order 2, where all elements on the plan $i=1$ are equal to zero, then based on definition \ref{defOrder3}:
\[
\det[B_{3 \times 3\times 3}]=
\det\begin{pmatrix} \left.
\begin{matrix}
\alpha \cdot a_{111} & \alpha \cdot a_{121} & \alpha \cdot a_{131}\\ 
a_{211} & a_{221} & a_{231}\\ 
a_{311} & a_{321} & a_{331}
\end{matrix}\right|
\begin{matrix}
\alpha \cdot a_{112} & \alpha \cdot a_{122} & \alpha \cdot a_{132}\\ 
a_{212} & a_{222} & a_{232}\\ 
a_{312} & a_{322} & a_{332}
\end{matrix}\left|
\begin{matrix}
\alpha \cdot a_{113} & \alpha \cdot a_{123} & \alpha \cdot a_{133}\\ 
a_{213} & a_{223} & a_{233}\\ 
a_{313} & a_{323} & a_{333}
\end{matrix}\right.
\end{pmatrix} 
\]
$$=\alpha \cdot a_{111} \cdot a_{222} \cdot a_{333} - \alpha \cdot a_{111} \cdot a_{232} \cdot a_{323} - \alpha \cdot a_{111} \cdot a_{223} \cdot a_{332} + \alpha \cdot a_{111} \cdot a_{233} \cdot a_{322}$$
$$-\alpha \cdot a_{112} \cdot a_{221} \cdot a_{333} + \alpha \cdot a_{112} \cdot a_{223} \cdot a_{331} + \alpha \cdot a_{112} \cdot a_{231} \cdot a_{323} - \alpha \cdot a_{112} \cdot a_{233} \cdot a_{321}$$
$$+\alpha \cdot a_{113} \cdot a_{221} \cdot a_{332} - \alpha \cdot a_{113} \cdot a_{222} \cdot a_{331} - \alpha \cdot a_{113} \cdot a_{231} \cdot a_{322} + \alpha \cdot a_{113} \cdot a_{232} \cdot a_{321}$$
$$-\alpha \cdot a_{121} \cdot a_{212} \cdot a_{333} + \alpha \cdot a_{121} \cdot a_{213} \cdot a_{332} + \alpha \cdot a_{121} \cdot a_{232} \cdot a_{313} - \alpha \cdot a_{121} \cdot a_{233} \cdot a_{312}$$
$$+\alpha \cdot a_{122} \cdot a_{211} \cdot a_{333} - \alpha \cdot a_{122} \cdot a_{213} \cdot a_{331} - \alpha \cdot a_{122} \cdot a_{231} \cdot a_{313} + \alpha \cdot a_{122} \cdot a_{233} \cdot a_{311}$$ 
$$-\alpha \cdot a_{123} \cdot a_{211} \cdot a_{332} + \alpha \cdot a_{123} \cdot a_{212} \cdot a_{331} + \alpha \cdot a_{123} \cdot a_{231} \cdot a_{312} - \alpha \cdot a_{123} \cdot a_{232} \cdot a_{311}$$
$$+\alpha \cdot a_{131} \cdot a_{212} \cdot a_{323} - \alpha \cdot a_{131} \cdot a_{213} \cdot a_{322} - \alpha \cdot a_{131} \cdot a_{222} \cdot a_{313} + \alpha \cdot a_{131} \cdot a_{223} \cdot a_{312}$$
$$-\alpha \cdot a_{132} \cdot a_{211} \cdot a_{323} + \alpha \cdot a_{132} \cdot a_{213} \cdot a_{321} + \alpha \cdot a_{132} \cdot a_{221} \cdot a_{313} - \alpha \cdot a_{132} \cdot a_{223} \cdot a_{311}$$
$$+\alpha \cdot a_{133} \cdot a_{211} \cdot a_{322} - \alpha \cdot a_{133} \cdot a_{212} \cdot a_{321} - \alpha \cdot a_{133} \cdot a_{221} \cdot a_{312} + \alpha \cdot a_{133} \cdot a_{222} \cdot a_{311} = \alpha \cdot \det[A_{3 \times 3\times 3}]$$

2. For plan $i=2$:

Let $A$ be cubic-matrix of order 2, where all elements on the plan $i=2$ are equal to zero, then based on definition \ref{defOrder3}:
\[
\det[B_{3 \times 3\times 3}]=
\det\begin{pmatrix} \left.
\begin{matrix}
a_{111} & a_{121} & a_{131}\\ 
\alpha \cdot a_{211} & \alpha \cdot a_{221} & \alpha \cdot a_{231}\\ 
a_{311} & a_{321} & a_{331}
\end{matrix}\right|
\begin{matrix}
a_{112} & a_{122} & a_{132}\\ 
\alpha \cdot a_{212} & \alpha \cdot a_{222} & \alpha \cdot a_{232}\\ 
a_{312} & a_{322} & a_{332}
\end{matrix}\left|
\begin{matrix}
a_{113} & a_{123} & a_{133}\\ 
\alpha \cdot a_{213} & \alpha \cdot a_{223} & \alpha \cdot a_{233}\\ 
a_{313} & a_{323} & a_{333}
\end{matrix}\right.
\end{pmatrix} 
\]
$$=a_{111} \cdot \alpha \cdot a_{222} \cdot a_{333} - a_{111} \cdot \alpha \cdot a_{232} \cdot a_{323} - a_{111} \cdot \alpha \cdot a_{223} \cdot a_{332} + a_{111} \cdot \alpha \cdot a_{233} \cdot a_{322}$$
$$-a_{112} \cdot \alpha \cdot a_{221} \cdot a_{333} + a_{112} \cdot \alpha \cdot a_{223} \cdot a_{331} + a_{112} \cdot \alpha \cdot a_{231} \cdot a_{323} - a_{112} \cdot \alpha \cdot a_{233} \cdot a_{321}$$
$$+a_{113} \cdot \alpha \cdot a_{221} \cdot a_{332} - a_{113} \cdot \alpha \cdot a_{222} \cdot a_{331} - a_{113} \cdot \alpha \cdot a_{231} \cdot a_{322} + a_{113} \cdot \alpha \cdot a_{232} \cdot a_{321}$$
$$-a_{121} \cdot \alpha \cdot a_{212} \cdot a_{333} + a_{121} \cdot \alpha \cdot a_{213} \cdot a_{332} + a_{121} \cdot \alpha \cdot a_{232} \cdot a_{313} - a_{121} \cdot \alpha \cdot a_{233} \cdot a_{312}$$
$$+a_{122} \cdot \alpha \cdot a_{211} \cdot a_{333} - a_{122} \cdot \alpha \cdot a_{213} \cdot a_{331} - a_{122} \cdot \alpha \cdot a_{231} \cdot a_{313} + a_{122} \cdot \alpha \cdot a_{233} \cdot a_{311}$$
$$-a_{123} \cdot \alpha \cdot a_{211} \cdot a_{332} + a_{123} \cdot \alpha \cdot a_{212} \cdot a_{331} + a_{123} \cdot \alpha \cdot a_{231} \cdot a_{312} - a_{123} \cdot \alpha \cdot a_{232} \cdot a_{311}$$
$$+a_{131} \cdot \alpha \cdot a_{212} \cdot a_{323} - a_{131} \cdot \alpha \cdot a_{213} \cdot a_{322} - a_{131} \cdot \alpha \cdot a_{222} \cdot a_{313} + a_{131} \cdot \alpha \cdot a_{223} \cdot a_{312}$$
$$-a_{132} \cdot \alpha \cdot a_{211} \cdot a_{323} + a_{132} \cdot \alpha \cdot a_{213} \cdot a_{321} + a_{132} \cdot \alpha \cdot a_{221} \cdot a_{313} - a_{132} \cdot \alpha \cdot a_{223} \cdot a_{311}$$
$$+a_{133} \cdot \alpha \cdot a_{211} \cdot a_{322} - a_{133} \cdot \alpha \cdot a_{212} \cdot a_{321} - a_{133} \cdot \alpha \cdot a_{221} \cdot a_{312} + a_{133} \cdot \alpha \cdot a_{222} \cdot a_{311} = \alpha \cdot \det[A_{3 \times 3\times 3}]$$

3. For plan $i=3$:

Let $A$ be cubic-matrix of order 2, where all elements on the plan $i=3$ are equal to zero, then based on definition \ref{defOrder3}:
\[
\det[B_{3 \times 3\times 3}]=
\det\begin{pmatrix} \left.
\begin{matrix}
a_{111} & a_{121} & a_{131}\\ 
a_{211} & a_{221} & a_{231}\\ 
\alpha \cdot a_{311} & \alpha \cdot a_{321} & \alpha \cdot a_{331}
\end{matrix}\right|
\begin{matrix}
a_{112} & a_{122} & a_{132}\\ 
a_{212} & a_{222} & a_{232}\\ 
\alpha \cdot a_{312} & \alpha \cdot a_{322} & \alpha \cdot a_{332}
\end{matrix}\left|
\begin{matrix}
a_{113} & a_{123} & a_{133}\\ 
a_{213} & a_{223} & a_{233}\\ 
\alpha \cdot a_{313} & \alpha \cdot a_{323} & \alpha \cdot a_{333}
\end{matrix}\right.
\end{pmatrix} 
\]
$$=a_{111} \cdot a_{222} \cdot \alpha \cdot a_{333} - a_{111} \cdot a_{232} \cdot \alpha \cdot a_{323} - a_{111} \cdot a_{223} \cdot \alpha \cdot a_{332} + a_{111} \cdot a_{233} \cdot \alpha \cdot a_{322}$$
$$-a_{112} \cdot a_{221} \cdot \alpha \cdot a_{333} + a_{112} \cdot a_{223} \cdot \alpha \cdot a_{331} + a_{112} \cdot a_{231} \cdot \alpha \cdot a_{323} - a_{112} \cdot a_{233} \cdot \alpha \cdot a_{321}$$
$$+a_{113} \cdot a_{221} \cdot \alpha \cdot a_{332} - a_{113} \cdot a_{222} \cdot \alpha \cdot a_{331} - a_{113} \cdot a_{231} \cdot \alpha \cdot a_{322} + a_{113} \cdot a_{232} \cdot \alpha \cdot a_{321}$$
$$-a_{121} \cdot a_{212} \cdot \alpha \cdot a_{333} + a_{121} \cdot a_{213} \cdot \alpha \cdot a_{332} + a_{121} \cdot a_{232} \cdot \alpha \cdot a_{313} - a_{121} \cdot a_{233} \cdot \alpha \cdot a_{312}$$
$$+a_{122} \cdot a_{211} \cdot \alpha \cdot a_{333} - a_{122} \cdot a_{213} \cdot \alpha \cdot a_{331} - a_{122} \cdot a_{231} \cdot \alpha \cdot a_{313} + a_{122} \cdot a_{233} \cdot \alpha \cdot a_{311}$$ 
$$-a_{123} \cdot a_{211} \cdot \alpha \cdot a_{332} + a_{123} \cdot a_{212} \cdot \alpha \cdot a_{331} + a_{123} \cdot a_{231} \cdot \alpha \cdot a_{312} - a_{123} \cdot a_{232} \cdot \alpha \cdot a_{311}$$
$$+a_{131} \cdot a_{212} \cdot \alpha \cdot a_{323} - a_{131} \cdot a_{213} \cdot \alpha \cdot a_{322} - a_{131} \cdot a_{222} \cdot \alpha \cdot a_{313} + a_{131} \cdot a_{223} \cdot \alpha \cdot a_{312}$$
$$-a_{132} \cdot a_{211} \cdot \alpha \cdot a_{323} + a_{132} \cdot a_{213} \cdot \alpha \cdot a_{321} + a_{132} \cdot a_{221} \cdot \alpha \cdot a_{313} - a_{132} \cdot a_{223} \cdot \alpha \cdot a_{311}$$
$$+a_{133} \cdot a_{211} \cdot \alpha \cdot a_{322} - a_{133} \cdot a_{212} \cdot \alpha \cdot a_{321} - a_{133} \cdot a_{221} \cdot \alpha \cdot a_{312} + a_{133} \cdot a_{222} \cdot \alpha \cdot a_{311} = \alpha \cdot \det[A_{3 \times 3\times 3}]$$

4. For plan $j=1$:

Let $A$ be cubic-matrix of order 2, where all elements on the plan $i=1$ are equal to zero, then based on definition \ref{defOrder3}:
\[
\det[B_{3 \times 3\times 3}]=
\det\begin{pmatrix} \left.
\begin{matrix}
\alpha \cdot a_{111} & a_{121} & a_{131}\\ 
\alpha \cdot a_{211} & a_{221} & a_{231}\\ 
\alpha \cdot a_{311} & a_{321} & a_{331}
\end{matrix}\right|
\begin{matrix}
\alpha \cdot a_{112} & a_{122} & a_{132}\\ 
\alpha \cdot a_{212} & a_{222} & a_{232}\\ 
\alpha \cdot a_{312} & a_{322} & a_{332}
\end{matrix}\left|
\begin{matrix}
\alpha \cdot a_{113} & a_{123} & a_{133}\\ 
\alpha \cdot a_{213} & a_{223} & a_{233}\\ 
\alpha \cdot a_{313} & a_{323} & a_{333}
\end{matrix}\right.
\end{pmatrix} 
\]
$$=\alpha \cdot a_{111} \cdot a_{222} \cdot a_{333} - \alpha \cdot a_{111} \cdot a_{232} \cdot a_{323} - \alpha \cdot a_{111} \cdot a_{223} \cdot a_{332} + \alpha \cdot a_{111} \cdot a_{233} \cdot a_{322}$$
$$-\alpha \cdot a_{112} \cdot a_{221} \cdot a_{333} + \alpha \cdot a_{112} \cdot a_{223} \cdot a_{331} + \alpha \cdot a_{112} \cdot a_{231} \cdot a_{323} - \alpha \cdot a_{112} \cdot a_{233} \cdot a_{321}$$
$$+\alpha \cdot a_{113} \cdot a_{221} \cdot a_{332} - \alpha \cdot a_{113} \cdot a_{222} \cdot a_{331} - \alpha \cdot a_{113} \cdot a_{231} \cdot a_{322} + \alpha \cdot a_{113} \cdot a_{232} \cdot a_{321}$$
$$-a_{121} \cdot \alpha \cdot a_{212} \cdot a_{333} + a_{121} \cdot \alpha \cdot a_{213} \cdot a_{332} + a_{121} \cdot a_{232} \cdot \alpha \cdot a_{313} - a_{121} \cdot a_{233} \cdot \alpha \cdot a_{312}$$
$$+a_{122} \cdot \alpha \cdot a_{211} \cdot a_{333} - a_{122} \cdot \alpha \cdot a_{213} \cdot a_{331} - a_{122} \cdot a_{231} \cdot \alpha \cdot a_{313} + a_{122} \cdot a_{233} \cdot \alpha \cdot a_{311}$$ 
$$-a_{123} \cdot \alpha \cdot a_{211} \cdot a_{332} + a_{123} \cdot \alpha \cdot a_{212} \cdot a_{331} + a_{123} \cdot a_{231} \cdot \alpha \cdot a_{312} - a_{123} \cdot a_{232} \cdot \alpha \cdot a_{311}$$
$$+a_{131} \cdot \alpha \cdot a_{212} \cdot a_{323} - a_{131} \cdot \alpha \cdot a_{213} \cdot a_{322} - a_{131} \cdot a_{222} \cdot \alpha \cdot a_{313} + a_{131} \cdot a_{223} \cdot \alpha \cdot a_{312}$$
$$-a_{132} \cdot \alpha \cdot a_{211} \cdot a_{323} + a_{132} \cdot \alpha \cdot a_{213} \cdot a_{321} + a_{132} \cdot a_{221} \cdot \alpha \cdot a_{313} - a_{132} \cdot a_{223} \cdot \alpha \cdot a_{311}$$
$$+a_{133} \cdot \alpha \cdot a_{211} \cdot a_{322} a_{133} \cdot \alpha \cdot a_{212} \cdot a_{321} - a_{133} \cdot a_{221} \cdot \alpha \cdot a_{312} + a_{133} \cdot a_{222} \cdot \alpha \cdot a_{311} = \alpha \cdot \det[A_{3 \times 3\times 3}]$$

5. For plan $j=2$:

Let $A$ be cubic-matrix of order 2, where all elements on the plan $i=2$ are equal to zero, then based on definition \ref{defOrder3}:
\[
\det[B_{3 \times 3\times 3}]=
\det\begin{pmatrix} \left.
\begin{matrix}
a_{111} & \alpha \cdot a_{121} & a_{131}\\ 
a_{211} & \alpha \cdot a_{221} & a_{231}\\ 
a_{311} & \alpha \cdot a_{321} & a_{331}
\end{matrix}\right|
\begin{matrix}
a_{112} & \alpha \cdot a_{122} & a_{132}\\ 
a_{212} & \alpha \cdot a_{222} & a_{232}\\ 
a_{312} & \alpha \cdot a_{322} & a_{332}
\end{matrix}\left|
\begin{matrix}
a_{113} & \alpha \cdot a_{123} & a_{133}\\ 
a_{213} & \alpha \cdot a_{223} & a_{233}\\ 
a_{313} & \alpha \cdot a_{323} & a_{333}
\end{matrix}\right.
\end{pmatrix} 
\]
$$=a_{111} \cdot \alpha \cdot a_{222} \cdot a_{333} - a_{111} \cdot a_{232} \cdot \alpha \cdot a_{323} - a_{111} \cdot \alpha \cdot a_{223} \cdot a_{332} + a_{111} \cdot a_{233} \cdot \alpha \cdot a_{322}$$
$$-a_{112} \cdot \alpha \cdot a_{221} \cdot a_{333} + a_{112} \cdot \alpha \cdot a_{223} \cdot a_{331} + a_{112} \cdot a_{231} \cdot \alpha \cdot a_{323} - a_{112} \cdot a_{233} \cdot \alpha \cdot a_{321}$$
$$+a_{113} \cdot \alpha \cdot a_{221} \cdot a_{332} - a_{113} \cdot \alpha \cdot a_{222} \cdot a_{331} - a_{113} \cdot a_{231} \cdot \alpha \cdot a_{322} + a_{113} \cdot a_{232} \cdot \alpha \cdot a_{321}$$
$$-\alpha \cdot a_{121} \cdot a_{212} \cdot a_{333} + \alpha \cdot a_{121} \cdot a_{213} \cdot a_{332} + \alpha \cdot a_{121} \cdot a_{232} \cdot a_{313} - \alpha \cdot a_{121} \cdot a_{233} \cdot a_{312}$$
$$+\alpha \cdot a_{122} \cdot a_{211} \cdot a_{333} - \alpha \cdot a_{122} \cdot a_{213} \cdot a_{331} - \alpha \cdot a_{122} \cdot a_{231} \cdot a_{313} + \alpha \cdot a_{122} \cdot a_{233} \cdot a_{311}$$ 
$$-\alpha \cdot a_{123} \cdot a_{211} \cdot a_{332} + \alpha \cdot a_{123} \cdot a_{212} \cdot a_{331} + \alpha \cdot a_{123} \cdot a_{231} \cdot a_{312} - \alpha \cdot a_{123} \cdot a_{232} \cdot a_{311}$$
$$+a_{131} \cdot a_{212} \cdot \alpha \cdot a_{323} - a_{131} \cdot a_{213} \cdot \alpha \cdot a_{322} - a_{131} \cdot \alpha \cdot a_{222} \cdot a_{313} + a_{131} \cdot \alpha \cdot a_{223} \cdot a_{312}$$
$$-a_{132} \cdot a_{211} \cdot \alpha \cdot a_{323} + a_{132} \cdot a_{213} \cdot \alpha \cdot a_{321} + a_{132} \cdot \alpha \cdot a_{221} \cdot a_{313} - a_{132} \cdot \alpha \cdot a_{223} \cdot a_{311}$$
$$+a_{133} \cdot a_{211} \cdot \alpha \cdot a_{322} - a_{133} \cdot a_{212} \cdot \alpha \cdot a_{321} - a_{133} \cdot \alpha \cdot a_{221} \cdot a_{312} + a_{133} \cdot \alpha \cdot a_{222} \cdot a_{311} = \alpha \cdot \det[A_{3 \times 3\times 3}]$$

6. For plan $j=3$:

Let $A$ be cubic-matrix of order 2, where all elements on the plan $i=3$ are equal to zero, then based on definition \ref{defOrder3}:
\[
\det[B_{3 \times 3\times 3}]=
\det\begin{pmatrix} \left.
\begin{matrix}
a_{111} & a_{121} & \alpha \cdot a_{131}\\ 
a_{211} & a_{221} & \alpha \cdot a_{231}\\ 
a_{311} & a_{321} & \alpha \cdot a_{331}
\end{matrix}\right|
\begin{matrix}
a_{112} & a_{122} & \alpha \cdot a_{132}\\ 
a_{212} & a_{222} & \alpha \cdot a_{232}\\ 
a_{312} & a_{322} & \alpha \cdot a_{332}
\end{matrix}\left|
\begin{matrix}
a_{113} & a_{123} & \alpha \cdot a_{133}\\ 
a_{213} & a_{223} & \alpha \cdot a_{233}\\ 
a_{313} & a_{323} & \alpha \cdot a_{333}
\end{matrix}\right.
\end{pmatrix} 
\]
$$=a_{111} \cdot a_{222} \cdot \alpha \cdot a_{333} - a_{111} \cdot \alpha \cdot a_{232} \cdot a_{323} - a_{111} \cdot a_{223} \cdot \alpha \cdot a_{332} + a_{111} \cdot \alpha \cdot a_{233} \cdot a_{322}$$
$$-a_{112} \cdot a_{221} \cdot \alpha \cdot a_{333} + a_{112} \cdot a_{223} \cdot \alpha \cdot a_{331} + a_{112} \cdot \alpha \cdot a_{231} \cdot a_{323} - a_{112} \cdot \alpha \cdot a_{233} \cdot a_{321}$$
$$+a_{113} \cdot a_{221} \cdot \alpha \cdot a_{332} - a_{113} \cdot a_{222} \cdot \alpha \cdot a_{331} - a_{113} \cdot \alpha \cdot a_{231} \cdot a_{322} + a_{113} \cdot \alpha \cdot a_{232} \cdot a_{321}$$
$$-a_{121} \cdot a_{212} \cdot \alpha \cdot a_{333} + a_{121} \cdot a_{213} \cdot \alpha \cdot a_{332} + a_{121} \cdot \alpha \cdot a_{232} \cdot a_{313} - a_{121} \cdot \alpha \cdot a_{233} \cdot a_{312}$$
$$+a_{122} \cdot a_{211} \cdot \alpha \cdot a_{333} - a_{122} \cdot a_{213} \cdot \alpha \cdot a_{331} - a_{122} \cdot \alpha \cdot a_{231} \cdot a_{313} + a_{122} \cdot \alpha \cdot a_{233} \cdot a_{311}$$ 
$$-a_{123} \cdot a_{211} \cdot \alpha \cdot a_{332} + a_{123} \cdot a_{212} \cdot \alpha \cdot a_{331} + a_{123} \cdot \alpha \cdot a_{231} \cdot a_{312} - a_{123} \cdot \alpha \cdot a_{232} \cdot a_{311}$$
$$+\alpha \cdot a_{131} \cdot a_{212} \cdot a_{323} - \alpha \cdot a_{131} \cdot a_{213} \cdot a_{322} - \alpha \cdot a_{131} \cdot a_{222} \cdot a_{313} + \alpha \cdot a_{131} \cdot a_{223} \cdot a_{312}$$
$$-\alpha \cdot a_{132} \cdot a_{211} \cdot a_{323} + \alpha \cdot a_{132} \cdot a_{213} \cdot a_{321} + \alpha \cdot a_{132} \cdot a_{221} \cdot a_{313} - \alpha \cdot a_{132} \cdot a_{223} \cdot a_{311}$$
$$+\alpha \cdot a_{133} \cdot a_{211} \cdot a_{322} - \alpha \cdot a_{133} \cdot a_{212} \cdot a_{321} - \alpha \cdot a_{133} \cdot a_{221} \cdot a_{312} + \alpha \cdot a_{133} \cdot a_{222} \cdot a_{311} = \alpha \cdot \det[A_{3 \times 3\times 3}]$$

7. For plan $k=1$:

Let $A$ be cubic-matrix of order 2, where all elements on the plan $k=1$ are equal to zero, then based on definition \ref{defOrder3}:
\[
\det[B_{3 \times 3\times 3}]=
\det\begin{pmatrix} \left.
\begin{matrix}
\alpha \cdot a_{111} & \alpha \cdot a_{121} & \alpha \cdot a_{131}\\ 
\alpha \cdot a_{211} & \alpha \cdot a_{221} & \alpha \cdot a_{231}\\ 
\alpha \cdot a_{311} & \alpha \cdot a_{321} & \alpha \cdot a_{331}
\end{matrix}\right|
\begin{matrix}
a_{112} & a_{122} & a_{132}\\ 
a_{212} & a_{222} & a_{232}\\ 
a_{312} & a_{322} & a_{332}
\end{matrix}\left|
\begin{matrix}
a_{113} & a_{123} & a_{133}\\ 
a_{213} & a_{223} & a_{233}\\ 
a_{313} & a_{323} & a_{333}
\end{matrix}\right.
\end{pmatrix} 
\]
$$=\alpha \cdot a_{111} \cdot a_{222} \cdot a_{333} - \alpha \cdot a_{111} \cdot a_{232} \cdot a_{323} - \alpha \cdot a_{111} \cdot a_{223} \cdot a_{332} + \alpha \cdot a_{111} \cdot a_{233} \cdot a_{322}$$
$$-a_{112} \cdot \alpha \cdot a_{221} \cdot a_{333} + a_{112} \cdot a_{223} \cdot \alpha \cdot a_{331} + a_{112} \cdot \alpha \cdot a_{231} \cdot a_{323} - a_{112} \cdot a_{233} \cdot \alpha \cdot a_{321}$$
$$+a_{113} \cdot \alpha \cdot a_{221} \cdot a_{332} - a_{113} \cdot a_{222} \cdot \alpha \cdot a_{331} - a_{113} \cdot \alpha \cdot a_{231} \cdot a_{322} + a_{113} \cdot a_{232} \cdot \alpha \cdot a_{321}$$
$$-\alpha \cdot a_{121} \cdot a_{212} \cdot a_{333} + \alpha \cdot a_{121} \cdot a_{213} \cdot a_{332} + \alpha \cdot a_{121} \cdot a_{232} \cdot a_{313} - \alpha \cdot a_{121} \cdot a_{233} \cdot a_{312}$$
$$+a_{122} \cdot \alpha \cdot a_{211} \cdot a_{333} - a_{122} \cdot a_{213} \cdot \alpha \cdot a_{331} - a_{122} \cdot \alpha \cdot a_{231} \cdot a_{313} + a_{122} \cdot a_{233} \cdot \alpha \cdot a_{311}$$ 
$$-a_{123} \cdot \alpha \cdot a_{211} \cdot a_{332} + a_{123} \cdot a_{212} \cdot \alpha \cdot a_{331} + a_{123} \cdot \alpha \cdot a_{231} \cdot a_{312} - a_{123} \cdot a_{232} \cdot \alpha \cdot a_{311}$$
$$+\alpha \cdot a_{131} \cdot a_{212} \cdot a_{323} - \alpha \cdot a_{131} \cdot a_{213} \cdot a_{322} - \alpha \cdot a_{131} \cdot a_{222} \cdot a_{313} + \alpha \cdot a_{131} \cdot a_{223} \cdot a_{312}$$
$$-a_{132} \cdot \alpha \cdot a_{211} \cdot a_{323} + a_{132} \cdot a_{213} \cdot \alpha \cdot a_{321} + a_{132} \cdot \alpha \cdot a_{221} \cdot a_{313} - a_{132} \cdot a_{223} \cdot \alpha \cdot a_{311}$$
$$ a_{133} \cdot \alpha \cdot a_{211} \cdot a_{322} - a_{133} \cdot a_{212} \cdot \alpha \cdot a_{321} - a_{133} \cdot \alpha \cdot a_{221} \cdot a_{312} + a_{133} \cdot a_{222} \cdot \alpha \cdot a_{311} = \alpha \cdot \det[A_{3 \times 3\times 3}]$$

8. For plan $k=2$:

Let $A$ be cubic-matrix of order 2, where all elements on the plan $k=2$ are equal to zero, then based on definition \ref{defOrder3}:
\[
\det[B_{3 \times 3\times 3}]=
\det\begin{pmatrix} \left.
\begin{matrix}
a_{111} & a_{121} & a_{131}\\ 
a_{211} & a_{221} & a_{231}\\ 
a_{311} & a_{321} & a_{331}
\end{matrix}\right|
\begin{matrix}
\alpha \cdot a_{112} & \alpha \cdot a_{122} & \alpha \cdot a_{132}\\ 
\alpha \cdot a_{212} & \alpha \cdot a_{222} & \alpha \cdot a_{232}\\ 
\alpha \cdot a_{312} & \alpha \cdot a_{322} & \alpha \cdot a_{332}
\end{matrix}\left|
\begin{matrix}
a_{113} & a_{123} & a_{133}\\ 
a_{213} & a_{223} & a_{233}\\ 
a_{313} & a_{323} & a_{333}
\end{matrix}\right.
\end{pmatrix} 
\]
$$=a_{111} \cdot \alpha \cdot a_{222} \cdot a_{333} - a_{111} \cdot \alpha \cdot a_{232} \cdot a_{323} - a_{111} \cdot a_{223} \cdot \alpha \cdot a_{332} + a_{111} \cdot a_{233} \cdot \alpha \cdot a_{322}$$
$$-\alpha \cdot a_{112} \cdot a_{221} \cdot a_{333} + \alpha \cdot a_{112} \cdot a_{223} \cdot a_{331} + \alpha \cdot a_{112} \cdot a_{231} \cdot a_{323} - \alpha \cdot a_{112} \cdot a_{233} \cdot a_{321}$$
$$+a_{113} \cdot a_{221} \cdot \alpha \cdot a_{332} - a_{113} \cdot \alpha \cdot a_{222} \cdot a_{331} - a_{113} \cdot a_{231} \cdot \alpha \cdot a_{322} + a_{113} \cdot \alpha \cdot a_{232} \cdot a_{321}$$
$$-a_{121} \cdot \alpha \cdot a_{212} \cdot a_{333} + a_{121} \cdot a_{213} \cdot \alpha \cdot a_{332} + a_{121} \cdot \alpha \cdot a_{232} \cdot a_{313} - a_{121} \cdot a_{233} \cdot \alpha \cdot a_{312}$$
$$+\alpha \cdot a_{122} \cdot a_{211} \cdot a_{333} - \alpha \cdot a_{122} \cdot a_{213} \cdot a_{331} - \alpha \cdot a_{122} \cdot a_{231} \cdot a_{313} + \alpha \cdot a_{122} \cdot a_{233} \cdot a_{311}$$ 
$$-a_{123} \cdot a_{211} \cdot \alpha \cdot a_{332} + a_{123} \cdot \alpha \cdot a_{212} \cdot a_{331} + a_{123} \cdot a_{231} \cdot \alpha \cdot a_{312} - a_{123} \cdot \alpha \cdot a_{232} \cdot a_{311}$$
$$+a_{131} \cdot \alpha \cdot a_{212} \cdot a_{323} - a_{131} \cdot a_{213} \cdot \alpha \cdot a_{322} - a_{131} \cdot \alpha \cdot a_{222} \cdot a_{313} + a_{131} \cdot a_{223} \cdot \alpha \cdot a_{312}$$
$$-\alpha \cdot a_{132} \cdot a_{211} \cdot a_{323} + \alpha \cdot a_{132} \cdot a_{213} \cdot a_{321} + \alpha \cdot a_{132} \cdot a_{221} \cdot a_{313} - \alpha \cdot a_{132} \cdot a_{223} \cdot a_{311}$$
$$+a_{133} \cdot a_{211} \cdot \alpha \cdot a_{322} - a_{133} \cdot \alpha \cdot a_{212} \cdot a_{321} - a_{133} \cdot a_{221} \cdot \alpha \cdot a_{312} + a_{133} \cdot \alpha \cdot a_{222} \cdot a_{311} = \alpha \cdot \det[A_{3 \times 3\times 3}]$$

9. For plan $k=3$:

Let $A$ be cubic-matrix of order 2, where all elements on the plan $k=3$ are equal to zero, then based on definition \ref{defOrder3}:
\[
\det[B_{3 \times 3\times 3}]=
\det\begin{pmatrix} \left.
\begin{matrix}
a_{111} & a_{121} & a_{131}\\ 
a_{211} & a_{221} & a_{231}\\ 
a_{311} & a_{321} & a_{331}
\end{matrix}\right|
\begin{matrix}
a_{112} & a_{122} & a_{132}\\ 
a_{212} & a_{222} & a_{232}\\ 
a_{312} & a_{322} & a_{332}
\end{matrix}\left|
\begin{matrix}
\alpha \cdot a_{113} & \alpha \cdot a_{123} & \alpha \cdot a_{133}\\ 
\alpha \cdot a_{213} & \alpha \cdot a_{223} & \alpha \cdot a_{233}\\ 
\alpha \cdot a_{313} & \alpha \cdot a_{323} & \alpha \cdot a_{333}
\end{matrix}\right.
\end{pmatrix} 
\]
$$=a_{111} \cdot a_{222} \cdot \alpha \cdot a_{333} - a_{111} \cdot a_{232} \cdot \alpha \cdot a_{323} - a_{111} \cdot \alpha \cdot a_{223} \cdot a_{332} + a_{111} \cdot \alpha \cdot a_{233} \cdot a_{322}$$
$$-a_{112} \cdot a_{221} \cdot \alpha \cdot a_{333} + a_{112} \cdot \alpha \cdot a_{223} \cdot a_{331} + a_{112} \cdot a_{231} \cdot \alpha \cdot a_{323} - a_{112} \cdot \alpha \cdot a_{233} \cdot a_{321}$$
$$+\alpha \cdot a_{113} \cdot a_{221} \cdot a_{332} - \alpha \cdot a_{113} \cdot a_{222} \cdot a_{331} - \alpha \cdot a_{113} \cdot a_{231} \cdot a_{322} + \alpha \cdot a_{113} \cdot a_{232} \cdot a_{321}$$
$$-a_{121} \cdot a_{212} \cdot \alpha \cdot a_{333} + a_{121} \cdot \alpha \cdot a_{213} \cdot a_{332} + a_{121} \cdot a_{232} \cdot \alpha \cdot a_{313} - a_{121} \cdot \alpha \cdot a_{233} \cdot a_{312}$$
$$+a_{122} \cdot a_{211} \cdot \alpha \cdot a_{333} - a_{122} \cdot \alpha \cdot a_{213} \cdot a_{331} - a_{122} \cdot a_{231} \cdot \alpha \cdot a_{313} + a_{122} \cdot \alpha \cdot a_{233} \cdot a_{311}$$ 
$$-\alpha \cdot a_{123} \cdot a_{211} \cdot a_{332} + \alpha \cdot a_{123} \cdot a_{212} \cdot a_{331} + \alpha \cdot a_{123} \cdot a_{231} \cdot a_{312} - \alpha \cdot a_{123} \cdot a_{232} \cdot a_{311}$$
$$+a_{131} \cdot a_{212} \cdot \alpha \cdot a_{323} - a_{131} \cdot \alpha \cdot a_{213} \cdot a_{322} - a_{131} \cdot a_{222} \cdot \alpha \cdot a_{313} + a_{131} \cdot \alpha \cdot a_{223} \cdot a_{312}$$
$$-a_{132} \cdot a_{211} \cdot \alpha \cdot a_{323} + a_{132} \cdot \alpha \cdot a_{213} \cdot a_{321} + a_{132} \cdot a_{221} \cdot \alpha \cdot a_{313} - a_{132} \cdot \alpha \cdot a_{223} \cdot a_{311}$$
$$+\alpha \cdot a_{133} \cdot a_{211} \cdot a_{322} - \alpha \cdot a_{133} \cdot a_{212} \cdot a_{321} - \alpha \cdot a_{133} \cdot a_{221} \cdot a_{312} + \alpha \cdot a_{133} \cdot a_{222} \cdot a_{311} = \alpha \cdot \det[A_{3 \times 3\times 3}]$$

\qed

\begin{example} Let's be $A$ a cubic matrix with order $2$ then we will obtain determinant a cubic-matrix $B$ from $A$ by multiplying any plan $i=1$ with scalar 3 and have,
\begin{equation*}
\det[B_{2 \times 2\times 2}]=\det
\begin{pmatrix}
\left.\begin{matrix}
3 \cdot 2 & 3 \cdot 1\\ 
3 & 5
\end{matrix}\right|
\begin{matrix}
3 \cdot 4 & 3 \cdot 7\\ 
3 & 2
\end{matrix}
\end{pmatrix}=3 \cdot 2\cdot 2 - 3 \cdot 4\cdot 5 - 3 \cdot 1\cdot 3 +3 \cdot 7\cdot 3=6.
\end{equation*}

If we compare with example 2, we can see that $|A|=3\cdot|B|$.
\end{example}

\begin{theorem} Let's be $A$ a cubic-matrix of order $2$ or order $3$, and $B$ be another cubic-matrix, which obtained from $A$ by interchanging the location of two consecutive "horizontal layer" in $i$-index, 
then $\det(A)=\det(B)$.
\end{theorem}

\proof 

Case 1. Let $B$ be cubic-matrix of order 2, where we have interchanged to "horizontal layer", then based on definition \ref{defOrder2}:

Let $B$ be cubic-matrix of order 2, where we have interchanged two "horizontal layer", then based on definition \ref{defOrder2}:
\[
\det[B_{2 \times 2\times 2}]=\det
\begin{pmatrix}
\left.\begin{matrix}
a_{211} & a_{221}\\ 
a_{111} & a_{121}
\end{matrix}\right|
\begin{matrix}
a_{212} & a_{222}\\ 
a_{112} & a_{122}
\end{matrix}
\end{pmatrix}\]
\[= a_{211}\cdot a_{122} - a_{212}\cdot a_{121} - a_{221}\cdot a_{112} + a_{222}\cdot a_{111} = \det[A_{2 \times 2\times 2}] 
\]

Case 2. Let $B$ be cubic-matrix of order 3, where we have interchanged to "horizontal layer", then based on definition \ref{defOrder2}:

1. Let $B$ be cubic-matrix of order 3, where we have interchanged two "horizontal layer" (first layer with the second layer), then based on definition \ref{defOrder3}:

\begin{equation*}
\det[B_{3 \times 3\times 3}]=
\det\begin{pmatrix} \left.
\begin{matrix}
a_{211} & a_{221} & a_{231}\\ 
a_{111} & a_{121} & a_{131}\\ 
a_{311} & a_{321} & a_{331}
\end{matrix}\right|
\begin{matrix}
a_{212} & a_{222} & a_{232}\\ 
a_{112} & a_{122} & a_{132}\\ 
a_{312} & a_{322} & a_{332}
\end{matrix}\left|
\begin{matrix}
a_{213} & a_{223} & a_{233}\\ 
a_{113} & a_{123} & a_{133}\\ 
a_{313} & a_{323} & a_{333}
\end{matrix}\right.
\end{pmatrix} 
\end{equation*}
$$=a_{211} \cdot a_{122} \cdot a_{333} - a_{211} \cdot a_{132} \cdot a_{323} - a_{211} \cdot a_{123} \cdot a_{332} + a_{211} \cdot a_{133} \cdot a_{322}$$
$$-a_{212} \cdot a_{121} \cdot a_{333} + a_{212} \cdot a_{123} \cdot a_{331} + a_{212} \cdot a_{131} \cdot a_{323} - a_{212} \cdot a_{133} \cdot a_{321}$$
$$+a_{213} \cdot a_{121} \cdot a_{332} - a_{213} \cdot a_{122} \cdot a_{331} - a_{213} \cdot a_{131} \cdot a_{322} + a_{213} \cdot a_{132} \cdot a_{321}$$
$$-a_{221} \cdot a_{112} \cdot a_{333} + a_{221} \cdot a_{113} \cdot a_{332} + a_{221} \cdot a_{132} \cdot a_{313} - a_{221} \cdot a_{133} \cdot a_{312}$$
$$+a_{222} \cdot a_{111} \cdot a_{333} - a_{222} \cdot a_{113} \cdot a_{331} - a_{222} \cdot a_{131} \cdot a_{313} + a_{222} \cdot a_{133} \cdot a_{311}$$ 
$$-a_{223} \cdot a_{111} \cdot a_{332} + a_{223} \cdot a_{112} \cdot a_{331} + a_{223} \cdot a_{131} \cdot a_{312} - a_{223} \cdot a_{132} \cdot a_{311}$$
$$+a_{231} \cdot a_{112} \cdot a_{323} - a_{231} \cdot a_{113} \cdot a_{322} - a_{231} \cdot a_{122} \cdot a_{313} + a_{231} \cdot a_{123} \cdot a_{312}$$
$$-a_{232} \cdot a_{111} \cdot a_{323} + a_{232} \cdot a_{113} \cdot a_{321} + a_{232} \cdot a_{121} \cdot a_{313} - a_{232} \cdot a_{123} \cdot a_{311}$$
$$+a_{233} \cdot a_{111} \cdot a_{322} - a_{233} \cdot a_{112} \cdot a_{311} - a_{233} \cdot a_{121} \cdot a_{312} + a_{233} \cdot a_{122} \cdot a_{311}$$

$$=a_{111} \cdot a_{222} \cdot a_{333} - a_{111} \cdot a_{232} \cdot a_{323} - a_{111} \cdot a_{223} \cdot a_{332} + a_{111} \cdot a_{233} \cdot a_{322}$$
$$-a_{112} \cdot a_{221} \cdot a_{333} + a_{112} \cdot a_{223} \cdot a_{331} + a_{112} \cdot a_{231} \cdot a_{323} - a_{112} \cdot a_{233} \cdot a_{321}$$
$$+a_{113} \cdot a_{221} \cdot a_{332} - a_{113} \cdot a_{222} \cdot a_{331} - a_{113} \cdot a_{231} \cdot a_{322} + a_{113} \cdot a_{232} \cdot a_{321}$$
$$-a_{121} \cdot a_{212} \cdot a_{333} + a_{121} \cdot a_{213} \cdot a_{332} + a_{121} \cdot a_{232} \cdot a_{313} - a_{121} \cdot a_{233} \cdot a_{312}$$
$$+a_{122} \cdot a_{211} \cdot a_{333} - a_{122} \cdot a_{213} \cdot a_{331} - a_{122} \cdot a_{231} \cdot a_{313} + a_{122} \cdot a_{233} \cdot a_{311}$$
$$-a_{123} \cdot a_{211} \cdot a_{332} + a_{123} \cdot a_{212} \cdot a_{331} + a_{123} \cdot a_{231} \cdot a_{312} - a_{123} \cdot a_{232} \cdot a_{311}$$
$$+a_{131} \cdot a_{212} \cdot a_{323} - a_{131} \cdot a_{213} \cdot a_{322} - a_{131} \cdot a_{222} \cdot a_{313} + a_{131} \cdot a_{223} \cdot a_{312}$$
$$-a_{132} \cdot a_{211} \cdot a_{323} + a_{132} \cdot a_{213} \cdot a_{321} + a_{132} \cdot a_{221} \cdot a_{313} - a_{132} \cdot a_{223} \cdot a_{311}$$
$$+a_{133} \cdot a_{211} \cdot a_{322} - a_{133} \cdot a_{212} \cdot a_{321} - a_{133} \cdot a_{221} \cdot a_{312} + a_{133} \cdot a_{222} \cdot a_{311} = \det[A_{3 \times 3\times 3}]$$

2. Let $B$ be cubic-matrix of order 3, where we have interchanged two "horizontal layer" (second layer with the third layer), then based on definition \ref{defOrder3}:

\begin{equation*}
\det[B_{3 \times 3\times 3}]=
\det\begin{pmatrix} \left.
\begin{matrix}
a_{111} & a_{121} & a_{131}\\ 
a_{311} & a_{321} & a_{331}\\ 
a_{211} & a_{221} & a_{231}
\end{matrix}\right|
\begin{matrix}
a_{112} & a_{122} & a_{132}\\ 
a_{312} & a_{322} & a_{332}\\ 
a_{212} & a_{222} & a_{232}
\end{matrix}\left|
\begin{matrix}
a_{113} & a_{123} & a_{133}\\ 
a_{313} & a_{323} & a_{333}\\ 
a_{213} & a_{223} & a_{233}
\end{matrix}\right.
\end{pmatrix} 
\end{equation*}
$$=a_{111} \cdot a_{322} \cdot a_{233} - a_{111} \cdot a_{332} \cdot a_{223} - a_{111} \cdot a_{323} \cdot a_{232} + a_{111} \cdot a_{333} \cdot a_{222}$$
$$-a_{112} \cdot a_{321} \cdot a_{233} + a_{112} \cdot a_{323} \cdot a_{231} + a_{112} \cdot a_{331} \cdot a_{223} - a_{112} \cdot a_{333} \cdot a_{221}$$
$$+a_{113} \cdot a_{321} \cdot a_{232} - a_{113} \cdot a_{322} \cdot a_{231} - a_{113} \cdot a_{331} \cdot a_{222} + a_{113} \cdot a_{332} \cdot a_{221}$$
$$-a_{121} \cdot a_{312} \cdot a_{233} + a_{121} \cdot a_{313} \cdot a_{232} + a_{121} \cdot a_{332} \cdot a_{213} - a_{121} \cdot a_{333} \cdot a_{212}$$
$$+a_{122} \cdot a_{311} \cdot a_{233} - a_{122} \cdot a_{313} \cdot a_{231} - a_{122} \cdot a_{331} \cdot a_{213} + a_{122} \cdot a_{333} \cdot a_{211}$$
$$-a_{123} \cdot a_{311} \cdot a_{232} + a_{123} \cdot a_{312} \cdot a_{231} + a_{123} \cdot a_{331} \cdot a_{212} - a_{123} \cdot a_{332} \cdot a_{211}$$
$$+a_{131} \cdot a_{312} \cdot a_{223} - a_{131} \cdot a_{313} \cdot a_{222} - a_{131} \cdot a_{322} \cdot a_{213} + a_{131} \cdot a_{323} \cdot a_{212}$$
$$-a_{132} \cdot a_{311} \cdot a_{223} + a_{132} \cdot a_{313} \cdot a_{221} + a_{132} \cdot a_{321} \cdot a_{213} - a_{132} \cdot a_{323} \cdot a_{211}$$
$$+a_{133} \cdot a_{311} \cdot a_{222} - a_{133} \cdot a_{312} \cdot a_{221} - a_{133} \cdot a_{321} \cdot a_{212} + a_{133} \cdot a_{322} \cdot a_{211}$$

$$=a_{111} \cdot a_{222} \cdot a_{333} - a_{111} \cdot a_{232} \cdot a_{323} - a_{111} \cdot a_{223} \cdot a_{332} + a_{111} \cdot a_{233} \cdot a_{322}$$
$$-a_{112} \cdot a_{221} \cdot a_{333} + a_{112} \cdot a_{223} \cdot a_{331} + a_{112} \cdot a_{231} \cdot a_{323} - a_{112} \cdot a_{233} \cdot a_{321}$$
$$+a_{113} \cdot a_{221} \cdot a_{332} - a_{113} \cdot a_{222} \cdot a_{331} - a_{113} \cdot a_{231} \cdot a_{322} + a_{113} \cdot a_{232} \cdot a_{321}$$
$$-a_{121} \cdot a_{212} \cdot a_{333} + a_{121} \cdot a_{213} \cdot a_{332} + a_{121} \cdot a_{232} \cdot a_{313} - a_{121} \cdot a_{233} \cdot a_{312}$$
$$+a_{122} \cdot a_{211} \cdot a_{333} - a_{122} \cdot a_{213} \cdot a_{331} - a_{122} \cdot a_{231} \cdot a_{313} + a_{122} \cdot a_{233} \cdot a_{311}$$
$$-a_{123} \cdot a_{211} \cdot a_{332} + a_{123} \cdot a_{212} \cdot a_{331} + a_{123} \cdot a_{231} \cdot a_{312} - a_{123} \cdot a_{232} \cdot a_{311}$$
$$+a_{131} \cdot a_{212} \cdot a_{323} - a_{131} \cdot a_{213} \cdot a_{322} - a_{131} \cdot a_{222} \cdot a_{313} + a_{131} \cdot a_{223} \cdot a_{312}$$
$$-a_{132} \cdot a_{211} \cdot a_{323} + a_{132} \cdot a_{213} \cdot a_{321} + a_{132} \cdot a_{221} \cdot a_{313} - a_{132} \cdot a_{223} \cdot a_{311}$$
$$+a_{133} \cdot a_{211} \cdot a_{322} - a_{133} \cdot a_{212} \cdot a_{321} - a_{133} \cdot a_{221} \cdot a_{312} + a_{133} \cdot a_{222} \cdot a_{311} = \det[A_{3 \times 3\times 3}]$$
\qed

\begin{example} Let $A$ be $2 \times 2 \times 2$ 3D determinant than we will obtain determinant $B$ from $A$ by interchanging location of two horizontal layer in $i$ index:

\begin{equation*}
\det[A_{2 \times 2\times 2}]=\det
\begin{pmatrix}
\left.\begin{matrix}
2 & 1\\ 
3 & 5
\end{matrix}\right|
\begin{matrix}
4 & 7\\ 
3 & 2
\end{matrix}
\end{pmatrix}=2\cdot 2 - 4\cdot 5 - 1\cdot 3 + 7\cdot 3=2.
\end{equation*}

Then:

\begin{equation*}
\det[A_{2 \times 2\times 2}]=\det
\begin{pmatrix}
\left.\begin{matrix}
3 & 5\\ 
2 & 1
\end{matrix}\right|
\begin{matrix}
3 & 2\\ 
4 & 7
\end{matrix}
\end{pmatrix}=3\cdot 7 - 3\cdot 1 - 5\cdot 4 + 2\cdot 2=2.
\end{equation*}

If we compare results with example 2, we can see that we have the same result. 
\end{example}

\begin{example} Let $A$ be a $3 \times 3 \times 3$ 3D determinant than we will obtain determinant $B$ from $A$ by interchanging location of two horizontal layer in $i$ index:

\begin{equation*}
\det[A_{3 \times 3\times 3}]=
\det\begin{pmatrix}
\left.
\begin{matrix}
1&4&2\\2&0&0\\0&4&2
\end{matrix}\right |
\begin{matrix}
3&1&3\\5&1&3\\3&2&0
\end{matrix}\left |
\begin{matrix}
2&1&0\\0&1&0\\2&1&0
\end{matrix} \right.
\end{pmatrix}=-3+6-4+24+24+10+10-4+6-6=63
\end{equation*}

Then:

\begin{equation*}
\det[B_{3 \times 3\times 3}]=
\det\begin{pmatrix}
\left.
\begin{matrix}
2&0&0\\1&4&2\\0&4&2
\end{matrix}\right |
\begin{matrix}
5&1&3\\3&1&3\\3&2&0
\end{matrix}\left |
\begin{matrix}
0&1&0\\2&1&0\\2&1&0
\end{matrix} \right.
\end{pmatrix}=6 -4 + 10 + 10 - 4 - 3 + 24 + 24 - 6 + 6 = 63
\end{equation*}

If we compare results with example 2, we can see that we have the same result.
\end{example}

\begin{theorem} Let's be $A$ a cubic-matrix of order $2$ or order $3$, and $B$ be another cubic-matrix, which obtained from $A$ by interchanging the location of two consecutive: "vertical page" in $j$-index or
"vertical layer" in $k$-index, 
then $\det(A)=-\det(B)$.
\end{theorem}

\proof 

Case 1. The cubic-matrix A of order 2, (and B has order 2), we will proof the case 1 for each "vertical page" and "vertical layer", as following:

1. For interchanging the location of two consecutive "vertical page":

Let $B$ be cubic-matrix of order 2, where we have interchanged two "vertical page", then based on definition \ref{defOrder2}:
\[
\det[B_{2 \times 2\times 2}]=\det
\begin{pmatrix}
\left.\begin{matrix}
a_{121} & a_{111}\\ 
a_{221} & a_{211}
\end{matrix}\right|
\begin{matrix}
a_{122} & a_{112}\\ 
a_{222} & a_{212}
\end{matrix}
\end{pmatrix}\]
\[=a_{121}\cdot a_{212} - a_{122}\cdot a_{211} - a_{111}\cdot a_{222} + a_{112}\cdot a_{221} = - \det[A_{2 \times 2\times 2}] 
\]

2. For interchanging the location of two consecutive "vertical layer":

Let $B$ be cubic-matrix of order 2, where we have interchanged two "vertical layer", then based on definition \ref{defOrder2}:
\[
\det[B_{2 \times 2\times 2}]=\det
\begin{pmatrix}
\left.\begin{matrix}
a_{112} & a_{122}\\ 
a_{212} & a_{222}
\end{matrix}\right|
\begin{matrix}
a_{111} & a_{121}\\ 
a_{211} & a_{221}
\end{matrix}
\end{pmatrix}\]
\[=a_{112}\cdot a_{221} - a_{111}\cdot a_{222} - a_{122}\cdot a_{211} + a_{121}\cdot a_{212} = - \det[A_{2 \times 2\times 2}] 
\]

Case 2. The cubic-matrix A of order 3, (and B has order 3), we will proof the case 2 for each "vertical page" and "vertical layer", as following:

1. For interchanging the location of two consecutive "vertical page":

Let $B$ be cubic-matrix of order 3, where we have interchanged two consecutive "vertical page" (First page with second page), then based on definition \ref{defOrder3}:

\[
\det[B_{3 \times 3\times 3}]=
\det\begin{pmatrix} \left.
\begin{matrix}
a_{121} & a_{111} & a_{131}\\ 
a_{221} & a_{211} & a_{231}\\ 
a_{321} & a_{311} & a_{331}
\end{matrix}\right|
\begin{matrix}
a_{122} & a_{112} & a_{132}\\ 
a_{222} & a_{212} & a_{232}\\ 
a_{322} & a_{312} & a_{332}
\end{matrix}\left|
\begin{matrix}
a_{123} & a_{113} & a_{133}\\ 
a_{223} & a_{213} & a_{233}\\ 
a_{323} & a_{313} & a_{333}
\end{matrix}\right.
\end{pmatrix} 
\]
$$=a_{121} \cdot a_{212} \cdot a_{333} - a_{121} \cdot a_{232} \cdot a_{313} - a_{121} \cdot a_{213} \cdot a_{332} + a_{121} \cdot a_{233} \cdot a_{312}$$
$$-a_{122} \cdot a_{211} \cdot a_{333} + a_{122} \cdot a_{213} \cdot a_{331} + a_{122} \cdot a_{231} \cdot a_{313} - a_{122} \cdot a_{233} \cdot a_{311}$$
$$+a_{123} \cdot a_{211} \cdot a_{332} - a_{123} \cdot a_{212} \cdot a_{331} - a_{123} \cdot a_{231} \cdot a_{312} + a_{123} \cdot a_{232} \cdot a_{311}$$
$$-a_{111} \cdot a_{222} \cdot a_{333} + a_{111} \cdot a_{223} \cdot a_{332} + a_{111} \cdot a_{232} \cdot a_{323} - a_{111} \cdot a_{233} \cdot a_{322}$$
$$+a_{112} \cdot a_{221} \cdot a_{333} - a_{112} \cdot a_{223} \cdot a_{331} - a_{112} \cdot a_{231} \cdot a_{323} + a_{112} \cdot a_{233} \cdot a_{321}$$
$$-a_{113} \cdot a_{221} \cdot a_{332} + a_{113} \cdot a_{222} \cdot a_{331} + a_{113} \cdot a_{231} \cdot a_{322} - a_{113} \cdot a_{232} \cdot a_{321}$$
$$+a_{131} \cdot a_{222} \cdot a_{313} - a_{131} \cdot a_{223} \cdot a_{312} - a_{131} \cdot a_{212} \cdot a_{323} + a_{131} \cdot a_{213} \cdot a_{322}$$
$$-a_{132} \cdot a_{221} \cdot a_{313} + a_{132} \cdot a_{223} \cdot a_{311} + a_{132} \cdot a_{211} \cdot a_{323} - a_{132} \cdot a_{213} \cdot a_{321}$$
$$+a_{133} \cdot a_{221} \cdot a_{312} - a_{133} \cdot a_{222} \cdot a_{311} - a_{133} \cdot a_{211} \cdot a_{322} + a_{133} \cdot a_{212} \cdot a_{321}$$

$$=-a_{111} \cdot a_{222} \cdot a_{333} + a_{111} \cdot a_{232} \cdot a_{323} + a_{111} \cdot a_{223} \cdot a_{332} - a_{111} \cdot a_{233} \cdot a_{322}$$
$$+a_{112} \cdot a_{221} \cdot a_{333} - a_{112} \cdot a_{223} \cdot a_{331} - a_{112} \cdot a_{231} \cdot a_{323} + a_{112} \cdot a_{233} \cdot a_{321}$$
$$-a_{113} \cdot a_{221} \cdot a_{332} + a_{113} \cdot a_{222} \cdot a_{331} + a_{113} \cdot a_{231} \cdot a_{322} - a_{113} \cdot a_{232} \cdot a_{321}$$
$$+a_{121} \cdot a_{212} \cdot a_{333} - a_{121} \cdot a_{213} \cdot a_{332} - a_{121} \cdot a_{232} \cdot a_{313} + a_{121} \cdot a_{233} \cdot a_{312}$$
$$-a_{122} \cdot a_{211} \cdot a_{333} + a_{122} \cdot a_{213} \cdot a_{331} + a_{122} \cdot a_{231} \cdot a_{313} - a_{122} \cdot a_{233} \cdot a_{311}$$
$$+a_{123} \cdot a_{211} \cdot a_{332} - a_{123} \cdot a_{212} \cdot a_{331} - a_{123} \cdot a_{231} \cdot a_{312} + a_{123} \cdot a_{232} \cdot a_{311}$$
$$-a_{131} \cdot a_{212} \cdot a_{323} + a_{131} \cdot a_{213} \cdot a_{322} + a_{131} \cdot a_{222} \cdot a_{313} - a_{131} \cdot a_{223} \cdot a_{312}$$
$$+a_{132} \cdot a_{211} \cdot a_{323} - a_{132} \cdot a_{213} \cdot a_{321} - a_{132} \cdot a_{221} \cdot a_{313} + a_{132} \cdot a_{223} \cdot a_{311}$$
$$-a_{133} \cdot a_{211} \cdot a_{322} + a_{133} \cdot a_{212} \cdot a_{321} + a_{133} \cdot a_{221} \cdot a_{312} - a_{133} \cdot a_{222} \cdot a_{311}$$

$$=-(a_{111} \cdot a_{222} \cdot a_{333} - a_{111} \cdot a_{232} \cdot a_{323} - a_{111} \cdot a_{223} \cdot a_{332} + a_{111} \cdot a_{233} \cdot a_{322}$$
$$-a_{112} \cdot a_{221} \cdot a_{333} + a_{112} \cdot a_{223} \cdot a_{331} + a_{112} \cdot a_{231} \cdot a_{323} - a_{112} \cdot a_{233} \cdot a_{321}$$
$$+a_{113} \cdot a_{221} \cdot a_{332} - a_{113} \cdot a_{222} \cdot a_{331} - a_{113} \cdot a_{231} \cdot a_{322} + a_{113} \cdot a_{232} \cdot a_{321}$$
$$-a_{121} \cdot a_{212} \cdot a_{333} + a_{121} \cdot a_{213} \cdot a_{332} + a_{121} \cdot a_{232} \cdot a_{313} - a_{121} \cdot a_{233} \cdot a_{312}$$
$$+a_{122} \cdot a_{211} \cdot a_{333} - a_{122} \cdot a_{213} \cdot a_{331} - a_{122} \cdot a_{231} \cdot a_{313} + a_{122} \cdot a_{233} \cdot a_{311}$$
$$-a_{123} \cdot a_{211} \cdot a_{332} + a_{123} \cdot a_{212} \cdot a_{331} + a_{123} \cdot a_{231} \cdot a_{312} - a_{123} \cdot a_{232} \cdot a_{311}$$
$$+a_{131} \cdot a_{212} \cdot a_{323} - a_{131} \cdot a_{213} \cdot a_{322} - a_{131} \cdot a_{222} \cdot a_{313} + a_{131} \cdot a_{223} \cdot a_{312}$$
$$-a_{132} \cdot a_{211} \cdot a_{323} + a_{132} \cdot a_{213} \cdot a_{321} + a_{132} \cdot a_{221} \cdot a_{313} - a_{132} \cdot a_{223} \cdot a_{311}$$
$$+a_{133} \cdot a_{211} \cdot a_{322} - a_{133} \cdot a_{212} \cdot a_{321} - a_{133} \cdot a_{221} \cdot a_{312} + a_{133} \cdot a_{222} \cdot a_{311}) = -\det[A_{3 \times 3\times 3}]$$

2. For interchanging the location of two consecutive "vertical page":

Let $B$ be cubic-matrix of order 3, where we have interchanged two consecutive "vertical page" (Second page with third page), then based on definition \ref{defOrder3}:

\[
\det[B_{3 \times 3\times 3}]=
\det\begin{pmatrix} \left.
\begin{matrix}
a_{111} & a_{131} & a_{121}\\ 
a_{211} & a_{231} & a_{221}\\ 
a_{311} & a_{331} & a_{321}
\end{matrix}\right|
\begin{matrix}
a_{112} & a_{132} & a_{122}\\ 
a_{212} & a_{232} & a_{222}\\ 
a_{312} & a_{332} & a_{322}
\end{matrix}\left|
\begin{matrix}
a_{113} & a_{133} & a_{123}\\ 
a_{213} & a_{233} & a_{223}\\ 
a_{313} & a_{333} & a_{323}
\end{matrix}\right.
\end{pmatrix} 
\]
$$=a_{111} \cdot a_{232} \cdot a_{323} - a_{111} \cdot a_{222} \cdot a_{333} - a_{111} \cdot a_{233} \cdot a_{322} + a_{111} \cdot a_{223} \cdot a_{332}$$
$$-a_{112} \cdot a_{231} \cdot a_{323} + a_{112} \cdot a_{233} \cdot a_{321} + a_{112} \cdot a_{221} \cdot a_{333} - a_{112} \cdot a_{223} \cdot a_{331}$$
$$+a_{113} \cdot a_{231} \cdot a_{322} - a_{113} \cdot a_{232} \cdot a_{321} - a_{113} \cdot a_{221} \cdot a_{332} + a_{113} \cdot a_{222} \cdot a_{331}$$
$$-a_{131} \cdot a_{212} \cdot a_{323} + a_{131} \cdot a_{213} \cdot a_{322} + a_{131} \cdot a_{222} \cdot a_{313} - a_{131} \cdot a_{223} \cdot a_{312}$$
$$+a_{132} \cdot a_{211} \cdot a_{323} - a_{132} \cdot a_{213} \cdot a_{321} - a_{132} \cdot a_{221} \cdot a_{313} + a_{132} \cdot a_{223} \cdot a_{311}$$
$$-a_{133} \cdot a_{211} \cdot a_{322} + a_{133} \cdot a_{212} \cdot a_{321} + a_{133} \cdot a_{221} \cdot a_{312} - a_{133} \cdot a_{222} \cdot a_{311}$$
$$+a_{121} \cdot a_{212} \cdot a_{333} - a_{121} \cdot a_{213} \cdot a_{332} - a_{121} \cdot a_{232} \cdot a_{313} + a_{121} \cdot a_{233} \cdot a_{312}$$
$$-a_{122} \cdot a_{211} \cdot a_{333} + a_{122} \cdot a_{213} \cdot a_{331} + a_{122} \cdot a_{231} \cdot a_{313} - a_{122} \cdot a_{233} \cdot a_{311}$$
$$+a_{123} \cdot a_{211} \cdot a_{332} - a_{123} \cdot a_{212} \cdot a_{331} - a_{123} \cdot a_{231} \cdot a_{312} + a_{123} \cdot a_{232} \cdot a_{311}$$

$$=-a_{111} \cdot a_{222} \cdot a_{333} + a_{111} \cdot a_{232} \cdot a_{323} + a_{111} \cdot a_{223} \cdot a_{332} - a_{111} \cdot a_{233} \cdot a_{322}$$
$$+a_{112} \cdot a_{221} \cdot a_{333} - a_{112} \cdot a_{223} \cdot a_{331} - a_{112} \cdot a_{231} \cdot a_{323} + a_{112} \cdot a_{233} \cdot a_{321}$$
$$-a_{113} \cdot a_{221} \cdot a_{332} + a_{113} \cdot a_{222} \cdot a_{331} + a_{113} \cdot a_{231} \cdot a_{322} - a_{113} \cdot a_{232} \cdot a_{321}$$
$$+a_{121} \cdot a_{212} \cdot a_{333} - a_{121} \cdot a_{213} \cdot a_{332} - a_{121} \cdot a_{232} \cdot a_{313} + a_{121} \cdot a_{233} \cdot a_{312}$$
$$-a_{122} \cdot a_{211} \cdot a_{333} + a_{122} \cdot a_{213} \cdot a_{331} + a_{122} \cdot a_{231} \cdot a_{313} - a_{122} \cdot a_{233} \cdot a_{311}$$
$$+a_{123} \cdot a_{211} \cdot a_{332} - a_{123} \cdot a_{212} \cdot a_{331} - a_{123} \cdot a_{231} \cdot a_{312} + a_{123} \cdot a_{232} \cdot a_{311}$$
$$-a_{131} \cdot a_{212} \cdot a_{323} + a_{131} \cdot a_{213} \cdot a_{322} + a_{131} \cdot a_{222} \cdot a_{313} - a_{131} \cdot a_{223} \cdot a_{312}$$
$$+a_{132} \cdot a_{211} \cdot a_{323} - a_{132} \cdot a_{213} \cdot a_{321} - a_{132} \cdot a_{221} \cdot a_{313} + a_{132} \cdot a_{223} \cdot a_{311}$$
$$-a_{133} \cdot a_{211} \cdot a_{322} + a_{133} \cdot a_{212} \cdot a_{321} + a_{133} \cdot a_{221} \cdot a_{312} - a_{133} \cdot a_{222} \cdot a_{311}$$

$$=-(a_{111} \cdot a_{222} \cdot a_{333} - a_{111} \cdot a_{232} \cdot a_{323} - a_{111} \cdot a_{223} \cdot a_{332} + a_{111} \cdot a_{233} \cdot a_{322}$$
$$-a_{112} \cdot a_{221} \cdot a_{333} + a_{112} \cdot a_{223} \cdot a_{331} + a_{112} \cdot a_{231} \cdot a_{323} - a_{112} \cdot a_{233} \cdot a_{321}$$
$$+a_{113} \cdot a_{221} \cdot a_{332} - a_{113} \cdot a_{222} \cdot a_{331} - a_{113} \cdot a_{231} \cdot a_{322} + a_{113} \cdot a_{232} \cdot a_{321}$$
$$-a_{121} \cdot a_{212} \cdot a_{333} + a_{121} \cdot a_{213} \cdot a_{332} + a_{121} \cdot a_{232} \cdot a_{313} - a_{121} \cdot a_{233} \cdot a_{312}$$
$$+a_{122} \cdot a_{211} \cdot a_{333} - a_{122} \cdot a_{213} \cdot a_{331} - a_{122} \cdot a_{231} \cdot a_{313} + a_{122} \cdot a_{233} \cdot a_{311}$$
$$-a_{123} \cdot a_{211} \cdot a_{332} + a_{123} \cdot a_{212} \cdot a_{331} + a_{123} \cdot a_{231} \cdot a_{312} - a_{123} \cdot a_{232} \cdot a_{311}$$
$$+a_{131} \cdot a_{212} \cdot a_{323} - a_{131} \cdot a_{213} \cdot a_{322} - a_{131} \cdot a_{222} \cdot a_{313} + a_{131} \cdot a_{223} \cdot a_{312}$$
$$-a_{132} \cdot a_{211} \cdot a_{323} + a_{132} \cdot a_{213} \cdot a_{321} + a_{132} \cdot a_{221} \cdot a_{313} - a_{132} \cdot a_{223} \cdot a_{311}$$
$$+a_{133} \cdot a_{211} \cdot a_{322} - a_{133} \cdot a_{212} \cdot a_{321} - a_{133} \cdot a_{221} \cdot a_{312} + a_{133} \cdot a_{222} \cdot a_{311}) = -\det[A_{3 \times 3\times 3}]$$

3. For interchanging the location of two consecutive "vertical layers":

Let $B$ be cubic-matrix of order 3, where we have interchanged two consecutive "vertical layers" (First layer with second layer), then based on definition \ref{defOrder3}:

\[
\det[B_{3 \times 3\times 3}]=
\det\begin{pmatrix} \left.
\begin{matrix}
a_{112} & a_{122} & a_{132}\\ 
a_{212} & a_{222} & a_{232}\\ 
a_{312} & a_{322} & a_{332}
\end{matrix}\right|
\begin{matrix}
a_{111} & a_{121} & a_{131}\\ 
a_{211} & a_{221} & a_{231}\\ 
a_{311} & a_{321} & a_{331}
\end{matrix}\left|
\begin{matrix}
a_{113} & a_{123} & a_{133}\\ 
a_{213} & a_{223} & a_{233}\\ 
a_{313} & a_{323} & a_{333}
\end{matrix}\right.
\end{pmatrix} 
\]
$$=a_{112} \cdot a_{221} \cdot a_{333} - a_{112} \cdot a_{231} \cdot a_{323} - a_{112} \cdot a_{223} \cdot a_{331} + a_{112} \cdot a_{233} \cdot a_{321}$$
$$-a_{111} \cdot a_{222} \cdot a_{333} + a_{111} \cdot a_{223} \cdot a_{332} + a_{111} \cdot a_{232} \cdot a_{323} - a_{111} \cdot a_{233} \cdot a_{322}$$
$$+a_{113} \cdot a_{222} \cdot a_{331} - a_{113} \cdot a_{221} \cdot a_{332} - a_{113} \cdot a_{232} \cdot a_{321} + a_{113} \cdot a_{231} \cdot a_{322}$$
$$-a_{122} \cdot a_{211} \cdot a_{333} + a_{122} \cdot a_{213} \cdot a_{331} + a_{122} \cdot a_{231} \cdot a_{313} - a_{122} \cdot a_{233} \cdot a_{311}$$
$$+a_{121} \cdot a_{212} \cdot a_{333} - a_{121} \cdot a_{213} \cdot a_{332} - a_{121} \cdot a_{232} \cdot a_{313} + a_{121} \cdot a_{233} \cdot a_{312}$$
$$-a_{123} \cdot a_{212} \cdot a_{331} + a_{123} \cdot a_{211} \cdot a_{332} + a_{123} \cdot a_{232} \cdot a_{311} - a_{123} \cdot a_{231} \cdot a_{312}$$
$$+a_{132} \cdot a_{211} \cdot a_{323} - a_{132} \cdot a_{213} \cdot a_{321} - a_{132} \cdot a_{221} \cdot a_{313} + a_{132} \cdot a_{223} \cdot a_{311}$$
$$-a_{131} \cdot a_{212} \cdot a_{323} + a_{131} \cdot a_{213} \cdot a_{322} + a_{131} \cdot a_{222} \cdot a_{313} - a_{131} \cdot a_{223} \cdot a_{312}$$
$$+a_{133} \cdot a_{212} \cdot a_{321} - a_{133} \cdot a_{211} \cdot a_{322} - a_{133} \cdot a_{222} \cdot a_{311} + a_{133} \cdot a_{221} \cdot a_{312}$$

$$=-a_{111} \cdot a_{222} \cdot a_{333} + a_{111} \cdot a_{232} \cdot a_{323} + a_{111} \cdot a_{223} \cdot a_{332} - a_{111} \cdot a_{233} \cdot a_{322}$$
$$+a_{112} \cdot a_{221} \cdot a_{333} - a_{112} \cdot a_{223} \cdot a_{331} - a_{112} \cdot a_{231} \cdot a_{323} + a_{112} \cdot a_{233} \cdot a_{321}$$
$$-a_{113} \cdot a_{221} \cdot a_{332} + a_{113} \cdot a_{222} \cdot a_{331} + a_{113} \cdot a_{231} \cdot a_{322} - a_{113} \cdot a_{232} \cdot a_{321}$$
$$+a_{121} \cdot a_{212} \cdot a_{333} - a_{121} \cdot a_{213} \cdot a_{332} - a_{121} \cdot a_{232} \cdot a_{313} + a_{121} \cdot a_{233} \cdot a_{312}$$
$$-a_{122} \cdot a_{211} \cdot a_{333} + a_{122} \cdot a_{213} \cdot a_{331} + a_{122} \cdot a_{231} \cdot a_{313} - a_{122} \cdot a_{233} \cdot a_{311}$$
$$+a_{123} \cdot a_{211} \cdot a_{332} - a_{123} \cdot a_{212} \cdot a_{331} - a_{123} \cdot a_{231} \cdot a_{312} + a_{123} \cdot a_{232} \cdot a_{311}$$
$$-a_{131} \cdot a_{212} \cdot a_{323} + a_{131} \cdot a_{213} \cdot a_{322} + a_{131} \cdot a_{222} \cdot a_{313} - a_{131} \cdot a_{223} \cdot a_{312}$$
$$+a_{132} \cdot a_{211} \cdot a_{323} - a_{132} \cdot a_{213} \cdot a_{321} - a_{132} \cdot a_{221} \cdot a_{313} + a_{132} \cdot a_{223} \cdot a_{311}$$
$$-a_{133} \cdot a_{211} \cdot a_{322} + a_{133} \cdot a_{212} \cdot a_{321} + a_{133} \cdot a_{221} \cdot a_{312} - a_{133} \cdot a_{222} \cdot a_{311}$$

$$=-(a_{111} \cdot a_{222} \cdot a_{333} - a_{111} \cdot a_{232} \cdot a_{323} - a_{111} \cdot a_{223} \cdot a_{332} + a_{111} \cdot a_{233} \cdot a_{322}$$
$$-a_{112} \cdot a_{221} \cdot a_{333} + a_{112} \cdot a_{223} \cdot a_{331} + a_{112} \cdot a_{231} \cdot a_{323} - a_{112} \cdot a_{233} \cdot a_{321}$$
$$+a_{113} \cdot a_{221} \cdot a_{332} - a_{113} \cdot a_{222} \cdot a_{331} - a_{113} \cdot a_{231} \cdot a_{322} + a_{113} \cdot a_{232} \cdot a_{321}$$
$$-a_{121} \cdot a_{212} \cdot a_{333} + a_{121} \cdot a_{213} \cdot a_{332} + a_{121} \cdot a_{232} \cdot a_{313} - a_{121} \cdot a_{233} \cdot a_{312}$$
$$+a_{122} \cdot a_{211} \cdot a_{333} - a_{122} \cdot a_{213} \cdot a_{331} - a_{122} \cdot a_{231} \cdot a_{313} + a_{122} \cdot a_{233} \cdot a_{311}$$
$$-a_{123} \cdot a_{211} \cdot a_{332} + a_{123} \cdot a_{212} \cdot a_{331} + a_{123} \cdot a_{231} \cdot a_{312} - a_{123} \cdot a_{232} \cdot a_{311}$$
$$+a_{131} \cdot a_{212} \cdot a_{323} - a_{131} \cdot a_{213} \cdot a_{322} - a_{131} \cdot a_{222} \cdot a_{313} + a_{131} \cdot a_{223} \cdot a_{312}$$
$$-a_{132} \cdot a_{211} \cdot a_{323} + a_{132} \cdot a_{213} \cdot a_{321} + a_{132} \cdot a_{221} \cdot a_{313} - a_{132} \cdot a_{223} \cdot a_{311}$$
$$+a_{133} \cdot a_{211} \cdot a_{322} - a_{133} \cdot a_{212} \cdot a_{321} - a_{133} \cdot a_{221} \cdot a_{312} + a_{133} \cdot a_{222} \cdot a_{311}) = -\det[A_{3 \times 3\times 3}]$$

4. For interchanging the location of two consecutive "vertical layers":

Let $B$ be cubic-matrix of order 3, where we have interchanged two consecutive "vertical layers" (Second layer with third layer), then based on definition \ref{defOrder3}:

\[
\det[B_{3 \times 3\times 3}]=
\det\begin{pmatrix} \left.
\begin{matrix}
a_{111} & a_{121} & a_{131}\\ 
a_{211} & a_{221} & a_{231}\\ 
a_{311} & a_{321} & a_{331}
\end{matrix}\right|
\begin{matrix}
a_{113} & a_{123} & a_{133}\\ 
a_{213} & a_{223} & a_{233}\\ 
a_{313} & a_{323} & a_{333}
\end{matrix}\left|
\begin{matrix}
a_{112} & a_{122} & a_{132}\\ 
a_{212} & a_{222} & a_{232}\\ 
a_{312} & a_{322} & a_{332}
\end{matrix}\right.
\end{pmatrix} 
\]
$$=a_{111} \cdot a_{223} \cdot a_{332} - a_{111} \cdot a_{233} \cdot a_{322} - a_{111} \cdot a_{222} \cdot a_{333} + a_{111} \cdot a_{232} \cdot a_{323}$$
$$-a_{113} \cdot a_{221} \cdot a_{332} + a_{113} \cdot a_{222} \cdot a_{331} + a_{113} \cdot a_{231} \cdot a_{322} - a_{113} \cdot a_{232} \cdot a_{321}$$
$$+a_{112} \cdot a_{221} \cdot a_{333} - a_{112} \cdot a_{223} \cdot a_{331} - a_{112} \cdot a_{231} \cdot a_{323} + a_{112} \cdot a_{233} \cdot a_{321}$$
$$-a_{121} \cdot a_{213} \cdot a_{332} + a_{121} \cdot a_{212} \cdot a_{333} + a_{121} \cdot a_{233} \cdot a_{312} - a_{121} \cdot a_{232} \cdot a_{313}$$
$$+a_{123} \cdot a_{211} \cdot a_{332} - a_{123} \cdot a_{212} \cdot a_{331} - a_{123} \cdot a_{231} \cdot a_{312} + a_{123} \cdot a_{232} \cdot a_{311}$$
$$-a_{122} \cdot a_{211} \cdot a_{333} + a_{122} \cdot a_{213} \cdot a_{331} + a_{122} \cdot a_{231} \cdot a_{313} - a_{122} \cdot a_{233} \cdot a_{311}$$
$$+a_{131} \cdot a_{213} \cdot a_{322} - a_{131} \cdot a_{212} \cdot a_{323} - a_{131} \cdot a_{223} \cdot a_{312} + a_{131} \cdot a_{222} \cdot a_{313}$$
$$-a_{133} \cdot a_{211} \cdot a_{322} + a_{133} \cdot a_{212} \cdot a_{321} + a_{133} \cdot a_{221} \cdot a_{312} - a_{133} \cdot a_{222} \cdot a_{311}$$
$$+a_{132} \cdot a_{211} \cdot a_{323} - a_{132} \cdot a_{213} \cdot a_{321} - a_{132} \cdot a_{221} \cdot a_{313} + a_{132} \cdot a_{223} \cdot a_{311}$$

$$=-a_{111} \cdot a_{222} \cdot a_{333} + a_{111} \cdot a_{232} \cdot a_{323} + a_{111} \cdot a_{223} \cdot a_{332} - a_{111} \cdot a_{233} \cdot a_{322}$$
$$+a_{112} \cdot a_{221} \cdot a_{333} - a_{112} \cdot a_{223} \cdot a_{331} - a_{112} \cdot a_{231} \cdot a_{323} + a_{112} \cdot a_{233} \cdot a_{321}$$
$$-a_{113} \cdot a_{221} \cdot a_{332} + a_{113} \cdot a_{222} \cdot a_{331} + a_{113} \cdot a_{231} \cdot a_{322} - a_{113} \cdot a_{232} \cdot a_{321}$$
$$+a_{121} \cdot a_{212} \cdot a_{333} - a_{121} \cdot a_{213} \cdot a_{332} - a_{121} \cdot a_{232} \cdot a_{313} + a_{121} \cdot a_{233} \cdot a_{312}$$
$$-a_{122} \cdot a_{211} \cdot a_{333} + a_{122} \cdot a_{213} \cdot a_{331} + a_{122} \cdot a_{231} \cdot a_{313} - a_{122} \cdot a_{233} \cdot a_{311}$$
$$+a_{123} \cdot a_{211} \cdot a_{332} - a_{123} \cdot a_{212} \cdot a_{331} - a_{123} \cdot a_{231} \cdot a_{312} + a_{123} \cdot a_{232} \cdot a_{311}$$
$$-a_{131} \cdot a_{212} \cdot a_{323} + a_{131} \cdot a_{213} \cdot a_{322} + a_{131} \cdot a_{222} \cdot a_{313} - a_{131} \cdot a_{223} \cdot a_{312}$$
$$+a_{132} \cdot a_{211} \cdot a_{323} - a_{132} \cdot a_{213} \cdot a_{321} - a_{132} \cdot a_{221} \cdot a_{313} + a_{132} \cdot a_{223} \cdot a_{311}$$
$$-a_{133} \cdot a_{211} \cdot a_{322} + a_{133} \cdot a_{212} \cdot a_{321} + a_{133} \cdot a_{221} \cdot a_{312} - a_{133} \cdot a_{222} \cdot a_{311}$$

$$=-(a_{111} \cdot a_{222} \cdot a_{333} - a_{111} \cdot a_{232} \cdot a_{323} - a_{111} \cdot a_{223} \cdot a_{332} + a_{111} \cdot a_{233} \cdot a_{322}$$
$$-a_{112} \cdot a_{221} \cdot a_{333} + a_{112} \cdot a_{223} \cdot a_{331} + a_{112} \cdot a_{231} \cdot a_{323} - a_{112} \cdot a_{233} \cdot a_{321}$$
$$+a_{113} \cdot a_{221} \cdot a_{332} - a_{113} \cdot a_{222} \cdot a_{331} - a_{113} \cdot a_{231} \cdot a_{322} + a_{113} \cdot a_{232} \cdot a_{321}$$
$$-a_{121} \cdot a_{212} \cdot a_{333} + a_{121} \cdot a_{213} \cdot a_{332} + a_{121} \cdot a_{232} \cdot a_{313} - a_{121} \cdot a_{233} \cdot a_{312}$$
$$+a_{122} \cdot a_{211} \cdot a_{333} - a_{122} \cdot a_{213} \cdot a_{331} - a_{122} \cdot a_{231} \cdot a_{313} + a_{122} \cdot a_{233} \cdot a_{311}$$
$$-a_{123} \cdot a_{211} \cdot a_{332} + a_{123} \cdot a_{212} \cdot a_{331} + a_{123} \cdot a_{231} \cdot a_{312} - a_{123} \cdot a_{232} \cdot a_{311}$$
$$+a_{131} \cdot a_{212} \cdot a_{323} - a_{131} \cdot a_{213} \cdot a_{322} - a_{131} \cdot a_{222} \cdot a_{313} + a_{131} \cdot a_{223} \cdot a_{312}$$
$$-a_{132} \cdot a_{211} \cdot a_{323} + a_{132} \cdot a_{213} \cdot a_{321} + a_{132} \cdot a_{221} \cdot a_{313} - a_{132} \cdot a_{223} \cdot a_{311}$$
$$+a_{133} \cdot a_{211} \cdot a_{322} - a_{133} \cdot a_{212} \cdot a_{321} - a_{133} \cdot a_{221} \cdot a_{312} + a_{133} \cdot a_{222} \cdot a_{311}) = -\det[A_{3 \times 3\times 3}]$$

\qed

\begin{example} Let $A$ be $2 \times 2 \times 2$ 3D determinant than we will obtain determinant $B$ from $A$ by interchanging location of two plans in $j$ index:

\begin{equation*}
\det[A_{2 \times 2\times 2}]=\det
\begin{pmatrix}
\left.\begin{matrix}
2 & 1\\ 
3 & 5
\end{matrix}\right|
\begin{matrix}
4 & 7\\ 
3 & 2
\end{matrix}
\end{pmatrix}=2\cdot 2 - 4\cdot 5 - 1\cdot 3 + 7\cdot 3=2.
\end{equation*}

Then:

\begin{equation*}
\det[B_{2 \times 2\times 2}]=\det
\begin{pmatrix}
\left.\begin{matrix}
1 & 2\\ 
5 & 3
\end{matrix}\right|
\begin{matrix}
7 & 4\\ 
2 & 3
\end{matrix}
\end{pmatrix}=1\cdot 3 - 7\cdot 3 - 2\cdot 2 + 4\cdot 5=-2.
\end{equation*}
\end{example}

\begin{example} Let $A$ be a $3 \times 3 \times 3$ 3D determinant than we will obtain determinant $B$ from $A$ by interchanging location of two plans in $j$ index:

\begin{equation*}
\det[A_{3 \times 3\times 3}]=
\det\begin{pmatrix}
\left.
\begin{matrix}
1&4&2\\2&0&0\\0&4&2
\end{matrix}\right |
\begin{matrix}
3&1&3\\5&1&3\\3&2&0
\end{matrix}\left |
\begin{matrix}
2&1&0\\0&1&0\\2&1&0
\end{matrix} \right.
\end{pmatrix}=-3+6-4+24+24+10+10-4+6-6=63
\end{equation*}

Then:

\begin{equation*}
\det[A_{3 \times 3\times 3}]=
\det\begin{pmatrix}
\left.
\begin{matrix}
4&1&2\\0&2&0\\4&0&2
\end{matrix}\right |
\begin{matrix}
1&3&3\\1&5&3\\2&3&0
\end{matrix}\left |
\begin{matrix}
1&2&0\\1&0&0\\1&2&0
\end{matrix} \right.
\end{pmatrix}=-24-10-6+3-6-24+4+4+6-10=-63
\end{equation*}

If we compare results with example 2, we can see that $|A|=-|B|$
\end{example}

\begin{example} Let $A$ be $2 \times 2 \times 2$ 3D determinant than we will obtain determinant $B$ from $A$ by interchanging location of two plans in $k$ index:

\begin{equation*}
\det[A_{2 \times 2\times 2}]=\det
\begin{pmatrix}
\left.\begin{matrix}
2 & 1\\ 
3 & 5
\end{matrix}\right|
\begin{matrix}
4 & 7\\ 
3 & 2
\end{matrix}
\end{pmatrix}=2\cdot 2 - 4\cdot 5 - 1\cdot 3 + 7\cdot 3=2.
\end{equation*}

Then:

\begin{equation*}
\det[B_{2 \times 2\times 2}]=\det
\begin{pmatrix}
\left.\begin{matrix}
4 & 7\\ 
3 & 2
\end{matrix}\right|
\begin{matrix}
2 & 1\\ 
3 & 5
\end{matrix}
\end{pmatrix}=4\cdot 5 - 2\cdot 2 - 7\cdot 3 + 1\cdot 3=-2.
\end{equation*}

If we compare results with example 2, we can see that $|A|=-|B|$.
\end{example}

\begin{example} Let $A$ be a $3 \times 3 \times 3$ 3D determinant than we will obtain determinant $B$ from $A$ by interchanging location of two plans in $k$ index:

\begin{equation*}
\det[A_{3 \times 3\times 3}]=
\det\begin{pmatrix}
\left.
\begin{matrix}
1&4&2\\2&0&0\\0&4&2
\end{matrix}\right |
\begin{matrix}
3&1&3\\5&1&3\\3&2&0
\end{matrix}\left |
\begin{matrix}
2&1&0\\0&1&0\\2&1&0
\end{matrix} \right.
\end{pmatrix}=-3+6-4+24+24+10+10-4+6-6=63
\end{equation*}

Then:

\begin{equation*}
\det[B_{3 \times 3\times 3}]=
\det\begin{pmatrix}
\left.
\begin{matrix}
3&1&3\\5&1&3\\3&2&0
\end{matrix}\right |
\begin{matrix}
1&4&2\\2&0&0\\0&4&2
\end{matrix}\left |
\begin{matrix}
2&1&0\\0&1&0\\2&1&0
\end{matrix} \right.
\end{pmatrix}=-6+3-24-24-10-6+6-10-4-4=-63
\end{equation*}

If we compare results with example 3, we can see that $|A|=-|B|$.
\end{example}

\section{Algorithms implementation of Determinants for cubic-matrix of order 2 and 3} %$ $\\

\subsection{Computer algorithm for calculating determinant of cubic matrices of order 2 and order 3} $ $\\

In the following we have presented the pseudocode of algorithm for calculating the determinant of cubic matrices of order 2 and order 3based on the Definition \ref{defOrder2} and Definition \ref{defOrder3}.

\quad

\noindent\makebox[\linewidth]{\rule{\textwidth}{1.3pt}}

\textbf{P 1:} Determinant calculation of cubic matrices of second and third order

\noindent\makebox[\linewidth]{\rule{\textwidth}{1.3pt}}

\textbf{Step 1:} Determine the order of determinant:

\hspace{1cm} $[m,n,o] = size(A);$

\textbf{Step 2:} Checking if 3D matrix is cubic:

\hspace{0.5cm} if $m ~\sim n; m \sim= o; n \sim= o;$

\hspace{1cm} disp('A is not square, cannot calculate the determinant')

\hspace{1cm} $d = 0$;

\hspace{1cm} return

\hspace{0.5cm} end

\textbf{Step 3:} Checking if 3D matrix is higher than the 3rd order:

\hspace{0.5cm} if $m > 3;$

\hspace{1cm} disp('A is higher than the third order, hence can not be calculated.')

\hspace{1cm} $d = 0$;

\hspace{1cm} return

\hspace{0.5cm} end

\textbf{Step 4:} Initialize $d=0;$

\textbf{Step 5:} Handling base case.

\hspace{0.5cm} if $m == 1$

\hspace{1cm} $d = A;$

\hspace{1cm} return

\hspace{0.5cm} end

\textbf{Step 6:} Check if A is of second order.

\hspace{0.5cm} $d = A(1,1,1) \ast A(2,2,2) - A(1,1,2) \ast A(2,2,1) - A(1,2,1) \ast A(2,1,2) + A(1,2,2) \ast A(2,1,1)$

\textbf{Step 7:} Check if A is of third order.

\hspace{0.5cm} $d=$

\hspace{0.5cm} $=A(1,1,1) \ast A(2,2,2) \ast A(3,3,3) - A(1,1,1) \ast A(2,3,2) \ast A(3,2,3) - A(1,1,1) \ast A(2,2,3) \ast A(3,3,2)$

\hspace{0.5cm} $+A(1,1,1) \ast A(2,3,3) \ast A(3,2,2) - A(1,1,2) \ast A(2,2,1) \ast A(3,3,3) + A(1,1,2) \ast A(2,2,3) \ast A(3,3,1)$

\hspace{0.5cm} $+A(1,1,2) \ast A(2,3,1) \ast A(3,2,3) - A(1,1,2) \ast A(2,3,3) \ast A(3,2,1) + A(1,1,3) \ast A(2,2,1) \ast A(3,3,2)$

\hspace{0.5cm} $-A(1,1,3) \ast A(2,2,2) \ast A(3,3,1) - A(1,1,3) \ast A(2,3,1) \ast A(3,2,2) + A(1,1,3) \ast A(2,3,2) \ast A(3,2,1)$

\hspace{0.5cm} $-A(1,2,1) \ast A(2,1,2) \ast A(3,3,3) + A(1,2,1) \ast A(2,1,3) \ast A(3,3,2) + A(1,2,1) \ast A(2,3,2) \ast A(3,1,3)$

\hspace{0.5cm} $-A(1,2,1) \ast A(2,3,3) \ast A(3,1,2) + A(1,2,2) \ast A(2,1,1) \ast A(3,3,3) - A(1,2,2) \ast A(2,1,3) \ast A(3,3,1)$

\hspace{0.5cm} $-A(1,2,2) \ast A(2,3,1) \ast A(3,1,3) + A(1,2,2) \ast A(2,3,3) \ast A(3,1,1) - A(1,2,3) \ast A(2,1,1) \ast A(3,3,2)$

\hspace{0.5cm} $+A(1,2,3) \ast A(2,1,2) \ast A(3,3,1) + A(1,2,3) \ast A(2,3,1) \ast A(3,1,2) - A(1,2,3) \ast A(2,3,2) \ast A(3,1,1)$

\hspace{0.5cm} $+A(1,3,1) \ast A(2,1,2) \ast A(3,2,3) - A(1,3,1) \ast A(2,1,3) \ast A(3,2,2) - A(1,3,1) \ast A(2,2,2) \ast A(3,1,3)$

\hspace{0.5cm} $+A(1,3,1) \ast A(2,2,3) \ast A(3,1,2) - A(1,3,2) \ast A(2,1,1) \ast A(3,2,3) + A(1,3,2) \ast A(2,1,3) \ast A(3,2,1)$

\hspace{0.5cm} $+A(1,3,2) \ast A(2,2,1) \ast A(3,1,3) - A(1,3,2) \ast A(2,2,3) \ast A(3,1,1) + A(1,3,3) \ast A(2,1,1) \ast A(3,2,2)$

\hspace{0.5cm} $-A(1,3,3) \ast A(2,1,2) \ast A(3,2,1) - A(1,3,3) \ast A(2,2,1) \ast A(3,1,2) + A(1,3,3) \ast A(2,2,2) \ast A(3,1,1)$

\textbf{Step 8:} Return the result of 3D determinant.

\noindent\makebox[\linewidth]{\rule{\textwidth}{1.3pt}}

\quad

\subsection{Optimized version of computer algorithm} $ $\\

The above-mention algorithm is hard-coded exactly as it is prescribet in the Dfinition \ref{defOrder2} and Definition \ref{defOrder3}, it can be optimized further with nested-loop. Hence, in the following we have presented another optimized version of above algorithm which gives the same result.

\quad

\noindent\makebox[\linewidth]{\rule{\textwidth}{1.3pt}}

\textbf{P 2:} Optimized algorithm for determinant calculation of cubic matrices of second and third order

\noindent\makebox[\linewidth]{\rule{\textwidth}{1.3pt}}

\textbf{Step 1:} Determine the order of determinant:

\hspace{1cm} $[m,n,o] = size(A);$

\textbf{Step 2:} Checking if 3D matrix is cubic:

\hspace{0.5cm} if $m ~\sim n; m \sim= o; n \sim= o;$

\hspace{1cm} disp('A is not square, cannot calculate the determinant')

\hspace{1cm} $d = 0$;

\hspace{1cm} return

\hspace{0.5cm} end

\textbf{Step 3:} Checking if 3D matrix is higher than the 3rd order:

\hspace{0.5cm} if $m > 3;$

\hspace{1cm} disp('A is higher than the third order, hence can not be calculated.')

\hspace{1cm} $d = 0$;

\hspace{1cm} return

\hspace{0.5cm} end

\textbf{Step 4:} Initialize $d=0;$

\textbf{Step 5:} Handling base case.

\hspace{0.5cm} if $m == 1$

\hspace{1cm} $d = A;$

\hspace{1cm} return

\hspace{0.5cm} end

\textbf{Step 6:} Check if A is of second order.

\hspace{0.5cm} Create loop for j from 1 to 2

\hspace{1cm} Create loop for k from 1 to 2

\hspace{1.5cm} Calculate determinant:

\hspace{1.8cm} $d = d + (-1)^{\wedge}(2+j+k) \ast A(1,j,k) \ast det\_3D( A(2,[1:j-1 j+1:2],[1:k-1 k+1:2]) );$

\hspace{1.5cm} end

\hspace{1cm} end

\textbf{Step 7:} Check if A is of third order.

\hspace{0.5cm} Create loop for j from 1 to 3

\hspace{1cm} Create loop for k from 1 to 3

\hspace{1.5cm} Calculate determinant:

\hspace{1.8cm} $d = d + (-1)^{\wedge}(2+j+k) \ast A(1,j,k) \ast det\_3D( A(2:3,[1:j-1 j+1:3],[1:k-1 k+1:3]) );$

\hspace{1.5cm} end

\hspace{1cm} end

\textbf{Step 8:} Return the result of 3D determinant.

\noindent\makebox[\linewidth]{\rule{\textwidth}{1.3pt}}

\section{Declarations} $ $\\
\textbf{Funding:} No Funding.\\
\textbf{Authors' contributions:} The contribution of the authors is equal.\\
\textbf{Data availability statements:} This manuscript does not report data.\\
\textbf{Conflict of Interest Statement:} There is no conflict of interest with any funder.


\begin{thebibliography}{99} % sort your citations by alphabetical order
\setlength{\itemsep}{1pt}
\setlength{\parskip}{0pt}
\setlength{\parsep}{0pt}

\bibitem{Salihu1} A. Salihu, H. Snopce, A. Luma and J. Ajdari, "Optimization of Dodgson's Condensation Method for Rectangular determinant Calculations", \textit{Advanced Mathematical Models and Applications}, vol. 7, no. 3, pp. 264-274, 2022. http://jomardpublishing.com/UploadFiles/Files/journals/AMMAV1N1/V7N3/Salihu\_et\_al.pdf.

\bibitem{PetersZakaDyckAM} Peters, J.F., Zaka, O. Dyck fundamental group on arcwise-connected polygon cycles. {\em Afr. Mat.}. \textbf{34}, 31 (2023), https://doi.org/10.1007/s13370-023-01067-3

\bibitem{ZakaDilauto} Zaka, O. Dilations of line in itself as the automorphism of the skew-field constructed over in the same line in Desargues affine plane. {\em Applied Mathematical Sciences}. \textbf{13}, 231-237 (2019)

\bibitem{ZakaFilipi2016}Zaka, O., Filipi, K. The transform of a line of Desargues affine plane in an additive group of its points. {\em Int. J. Of Current Research}. \textbf{8}, 34983-34990 (2016)

\bibitem{FilipiZakaJusufi}Filipi, K., Zaka, O., Jusufi, A. The construction of a corp in the set of points in a line of Desargues affine plane. {\em Matematicki Bilten}. \textbf{43}, 1-23 (2019), ISSN 0351-336X (print), ISSN 1857–9914 (online)

\bibitem{ZakaCollineations} Zaka, O. A description of collineations-groups of an affine plane. {\em Libertas Mathematica (N.S.)}. \textbf{37}, 81-96 (2017), ISSN print: 0278 – 5307, ISSN online: 2182 – 567X, MR3828328

\bibitem{ZakaVertex} Zaka, O. Three Vertex and Parallelograms in the Affine Plane: Similarity and Addition Abelian Groups of Similarly n-Vertexes in the Desargues Affine Plane. {\em Mathematical Modelling And Applications}. \textbf{3}, 9-15 (2018), http://doi:10.11648/j.mma.20180301.12

\bibitem{ZakaThesisPhd} Zaka, O. Contribution to Reports of Some Algebraic Structures with Affine Plane Geometry and Applications. (Polytechnic University of Tirana,Tirana, Albania,2016), supervisor: K. Filipi, vii+113pp.

\bibitem{ZakaPetersIso} Orgest Zaka and James F. Peters. Isomorphic-dilations of the skew-fields constructed over parallel lines in the Desargues affine plane. {\em Balkan J. Geom. Appl.}. \textbf{25}, 141-157 (2020), www.mathem.pub.ro/bjga/v25n1/B25-1zk-ZBG89.pdf

\bibitem{ZakaPetersOrder} Orgest Zaka and James Francis Peters. Ordered line and skew-fields in the Desargues affine plane. {\em Balkan J. Geom. Appl.}. \textbf{26}, 141-156 (2021), www.mathem.pub.ro/bjga/v26n1/B26-1zb-ZBP43.pdf

\bibitem{ZakaMohammedSF} O. Zaka and M. A. Mohammed, "Skew-field of trace-preserving endomorphisms, of translation group in affine plane", \textit{Proyecciones (Antofagasta, On line)}, vol. 39, no. 4, pp. 823-850, Jul. 2020. https://doi.org/10.22199/issn.0717-6279-2020-04-0052

\bibitem{ZakaMohammedEndo} O. Zaka and M. A. Mohammed, "The endomorphisms algebra of translations group and associative unitary ring of trace-preserving endomorphisms in affine plane", \textit{Proyecciones (Antofagasta, On line)}, vol. 39, no. 4, pp. 821-834, Jul. 2020. https://doi.org/10.22199/issn.0717-6279-2020-04-0051

\bibitem{Salihu2} A. Salihu, H. Snopce, A. Luma and J. Ajdari, "Comparison of time complexity growth for different methods/algorithms for rectangular determinant calculations", \textit{ICRTEC 2023 - Proceedings: IEEE International Conference on Recent Trends in Electronics and Communication: Upcoming Technologies for Smart Systems}. https://doi.org/10.1109/ICRTEC56977.2023.10111874.

\bibitem{Salihu3} A. Salihu, H. Snopce, J. Ajdari and A. Luma, "Generalization of Dodgson’s condensation method for calculating determinant of rectangular matrices", \textit{International Conference on Electrical, Computer and Energy Technologies (ICECET)}. https://doi.org/10.1109/ICECET55527.2022.9873054.

\bibitem{Salihu4} A. Salihu, H. Snopce, A. Luma and J. Ajdari, "Time Complexity Analysis for Cullis/Radic and Dodgson’s Generalized/Modified Method for Rectangular Determinants Calculations", \textit{International Journal of Computers and Their Applications}, vol. 29, no. 4, pp. 236-246, 2022. http://isca-hq.org/Documents/Journal/Archive/2022/2022volume2904/2022volume290403.pdf.

\bibitem{Salihu5} A. Salihu and F. Marevci, "Chio's-like Method for Calculating the Rectangular (non-square) Determinants: Computer Algorithm Interpretation and Comparison", \textit{European Journal of Pure and Applied Mathematics}, vol. 14, no. 2, pp. 431-450, 2021. https://doi.org/10.29020/nybg.ejpam.v14i2.3920.

\bibitem{Salihu6} A. Salihu and F. Marevci, "Determinants Order Decrease/Increase for k Orders, Interpretation with Computer Algorithms and Comparison", \textit{International Journal of Mathematics and Computer Science}, vol. 14, no. 2, pp. 501-518, 2021. http://ijmcs.future-in-tech.net/14.2/R-Marecvi-Salihu.pdf.

\bibitem{Salihu7} A. Salihu, A. Jusufi and F. Salihu, "Comparison of Computer Execution Time of Cornice Determinant Calculation", \textit{International Journal of Mathematics and Computer Science}, vol. 14, pp. 9-16, 2019. http://ijmcs.future-in-tech.net/14.1/R-Salihu2.pdf.

\bibitem{Salihu8} A. Salihu, "A modern modification of Gjonbalaj-Salihu cornice determinant, transformation to semi-diagonal determinant", \textit{International Journal of Mathematics and Computer Science}, vol. 13, pp. 1330138, 2018. http://ijmcs.future-in-tech.net/13.2/R-Salihu.pdf.

\bibitem{zaka3DmatrixRing} ZAKA, O. (2017) 3D Matrix Ring with a “Common” Multiplication. Open Access Library Journal, 4, 1-11. doi: http://dx.doi.org/10.4236/oalib.1103593.

\bibitem{zaka3DGLnnp} {Zaka, Orgest}, {The general linear group of degree {{\(n\)}} for 3D matrices {{\(\mathrm{GL}(n;n;p;F)\)}}}. Libertas Mathematica, New Series. Lib. Math. (N.S.) 39, No. 1, 13--30 (2019; Zbl 1451.15007)

\bibitem{ArtinM} Artin, M. (1991) Algebra. Prentice Hall, Upper Saddle River.

\bibitem{BretscherO} Bretscher, O. (2005) Linear Algebra with Applications. 3rd Edition, Prentice Hall, Upper Saddle River

\bibitem{Schneide-Barker} Schneide, H. and Barker, G.P. (1973) Matrices and Linear Algebra (Dover Books on Mathematics). 2nd Revised Edition.

\bibitem{Lang} Lang, S. (1987) Linear Algebra. Springer-Verlag, Berlin, New York.

\bibitem{Amiri-etal} Amiri, M., Fathy, M., Bayat, M., Generalization of some determinantal identities for non-square matrices based on Radic's definition, TWMS J. Pure Appl. Math. 1, no. 2 (2010), 163–175.

\bibitem{Radic1} Radić, M., A definition of determinant of rectangular matrix, Glas. Mat. Ser. III 1(21) (1966), 17–22.

\bibitem{Radic2} Radić, M., About a determinant of rectangular 2 × n matrix and its geometric interpretation, Beitr¨age Algebra Geom. 46, no. 2 (2005), 321–349

\bibitem{MAKAREWICZetal} Anna Makarewicz, Piotr Pikuta, and Dominik Szalkowski. "Properties of the determinant of a rectangular matrix." Annales Universitatis Mariae Curie-Skłodowska, sectio A – Mathematica 68.1 (2014): null. <http://eudml.org/doc/289812>.


\bibitem{Milne-Thomson} Milne-Thomson, L. (1941). Determinant Expansions. The Mathematical Gazette, 25(265), 130-135. doi:10.2307/3607371

% ...



\end{thebibliography}
\end{document}